\documentclass[12pt]{article}
\usepackage{amsmath,amssymb,bm,color}

\newtheorem{thm}{Theorem}[subsection]
\newtheorem{lem}[thm]{Lemma}
\newtheorem{prop}[thm]{Proposition}
\newtheorem{cor}[thm]{Corollary}
\newtheorem{define}[thm]{Definition}
\newtheorem{ex}[thm]{Example}
\newtheorem{rem}[thm]{Remark}
\newtheorem{conj}[thm]{Conjecture}

\newcommand{\qed}{\hfill$\Box$\par}

\begin{document}
\baselineskip=20pt
\begin{center}
\textbf{\Large Path Model for Representations of Generalized Kac--Moody Algebras}
\end{center}
\begin{center}
{\large Motohiro Ishii}\\
Graduate School of Pure and Applied Sciences, University of Tsukuba, \\
Tsukuba, Ibaraki 305-8571, Japan \\
(e-mail: \verb/ishii731@math.tsukuba.ac.jp/)
\end{center}
\begin{quote}
\begin{center}
\textbf{Abstract}
\end{center}
In [\textbf{JL}], Joseph and Lamprou generalized 
Littelmann's path model for Kac--Moody algebras to the case of generalized Kac--Moody algebras.
We show that Joseph--Lamprou's path model can be embedded into 
Littelmann's path model for a certain Kac--Moody algebra constructed from the 
Borcherds--Cartan datum of a given generalized Kac--Moody algebra; note that
this is not an embedding of crystals. Using this embedding, we give a new proof of the 
isomorphism theorem for path crystals, obtained in [\textbf{JL}, \S 7.4]. 
Moreover, for Joseph--Lamprou's path crystals, we give a decomposition rule for tensor product
and a branching rule for restriction to Levi subalgebras. 
Also, we obtain a characterization of standard paths in terms of a certain monoid
which can be thought of as a generalization of a Weyl group. 
\end{quote}

\section{Introduction}
The quantized universal enveloping algebra 
$U_q (\mathfrak{g})$ associated with a symmetrizable Kac--Moody algebra $\mathfrak{g}$
was introduced independently by V. G. Drinfel'd ([\textbf{D}]) and M. Jimbo ([\textbf{Ji}]). 
It is a non-commutative and non-cocommutative Hopf algebra that can be thought of as 
a $q$-deformation of the universal enveloping algebra $U(\mathfrak{g})$ of $\mathfrak{g}$. 
Also, G. Lusztig ([\textbf{Lu}]) and M. Kashiwara ([\textbf{Kas1-4}])
independently developed the theory of crystal bases (or canonical bases).
In particular, crystal bases can be regarded as bases at the limit ``$q=0$", 
where many of the problems in representation theory
can be reduced to combinatorial ones. 

In [\textbf{JKK, JKKS, Kan}], a theory parallel to the above was developed
for generalized Kac--Moody algebras, which is a new class of infinite-dimensional Lie algebras
introduced by R. Borcherds in his study of the Conway--Norton monstrous moonshine ([\textbf{B1-2, CN}]).
The representation theory of generalized Kac--Moody algebras is very similar to that of Kac--Moody algebras, 
and many results for Kac--Moody algebras can be extended to the case of generalized Kac--Moody algebras.
However, there indeed exist some differences between these Lie algebras. 
In particular, generalized Kac--Moody algebras have \textit{imaginary simple roots}, 
which have norms at most $0$; they cause some pathological phenomena.

Motivated by the crystal basis theory, P. Littelmann introduced a 
path model for integrable highest weight modules over the quantized universal enveloping 
algebra of a symmetrizable Kac--Moody algebra ([\textbf{Li1-2}]).
Soon after, it was shown in [\textbf{Jo, Kas5}] that Littelmann's path model provides a combinatorial realization 
of crystal bases of integrable highest weight modules. In Littelmann's path model, 
each crystal basis element is realized as a
\textit{Lakshmibai--Seshadri path} (\textit{LS path} for short), which is a pair $(\bm{w}; \bm{a})$ of a 
decreasing sequence of Weyl group elements 
$\bm{w}=(w_1 \ge w_2 \ge \cdots \ge w_k )$ (in the Bruhat order) 
and an increasing sequence of rational numbers 
$\bm{a}=(0=a_0 < a_1 < \cdots < a_k =1)$,
with certain integrality conditions. In [\textbf{JL}], A. Joseph and P. Lamprou generalized Littelmann's path model 
to the case of generalized Kac--Moody algebras.
They introduced \textit{generalized Lakshmibai--Seshadri paths} (\textit{GLS paths} for short), 
and considered the crystal $\mathbb{B}(\lambda )$ of all GLS paths 
of shape $\lambda $. Moreover, they showed that the (formal) character of $\mathbb{B}(\lambda )$
coincides with the generating function for the weights given by the Weyl--Kac--Borcherds 
character formula for the integrable highest weight
module with highest weight $\lambda $ over a generalized Kac--Moody algebra.

In this paper, we introduce a new tool, which we call an \textit{embedding of path models}. 
This is an embedding of the path crystal $\mathbb{B}(\lambda )$ of GLS paths of shape $\lambda $ for  
a given generalized Kac--Moody algebra $\mathfrak{g}$ into the path crystal 
$\widetilde{\mathbb{B}}(\Tilde{\lambda })$ of ordinary LS paths of shape $\Tilde{\lambda }$ for a certain 
Kac--Moody algebra $\Tilde{\mathfrak{g}}$ (Proposition \ref{4.1.2}); 
note that this is not an embedding of crystals.
As an application of this embedding, we give a new proof of the isomorphism theorem for path crystals, 
obtained in [\textbf{JL}, \S 7.4]; in fact, this is a special case of the following more general result:
\begin{thm}
Let $\lambda \in P^+$ be a dominant integral weight, and let
$\lambda _1 ,\ldots ,\lambda _n$ be integral weights 
in the Tits cone $\mathcal{W}P^+$ (see \S 4) 
such that $\lambda =\lambda _1 +\cdots +\lambda _n$, and such that
the concatenation $\pi _{\lambda _1}\otimes \cdots 
\otimes \pi _{\lambda _n }$ of the straight-line paths $\pi _{\lambda _1},\ldots ,\pi _{\lambda _n}$
is a dominant path, i.e., $\bigl( \pi _{\lambda _1}\otimes \cdots \otimes \pi _{\lambda _n }\bigl) (t) 
\in \sum _{\eta \in P^+} \mathbb{R}_{\ge 0}\eta $ for all $t\in [0,1]$.
Then, we have an isomorphism of crystals:
\begin{center}
$\mathbb{B}(\lambda ) \cong \mathcal{F}
(\pi _{\lambda _1}\otimes \cdots \otimes \pi _{\lambda _n }),$
\end{center}
where $\mathcal{F}$ is the monoid generated by root operators $f_i ,\ i\in I$ (see \S 3.1).
\label{1.0.1}
\end{thm}

First, we associate a certain monoid $\mathcal{W}$
with a given Borcherds--Cartan matrix, which can be 
thought of as a generalization of a Weyl group; indeed, this monoid has many nice combinatorial 
properties such as the (Strong) \textit{Exchange Property}, \textit{Deletion Property}, 
\textit{Word Property}, and \textit{Lifting Property}.
Then, we give a characterization of standard paths in terms of the \textit{Bruhat order} on the 
monoid $\mathcal{W}$ (Theorem \ref{5.2.2}), generalizing the corresponding result [\textbf{Li3}, Theorem 10.1]
for Kac--Moody algebras, due to P. Littelmann. By using this characterization of standard paths, 
we recover the crystal isomorphism theorem between Joseph--Lamprou's path crystal and the crystal basis 
of the integrable highest weight $U_q (\mathfrak{g})$-module, which is 
a special case of Theorem \ref{1.0.1}, as follows.
\begin{thm}\ 
For each dominant integral weight $\lambda \in P^+$, 
we have the following isomorphism of crystals:
\begin{center}
$\mathbb{B}(\lambda ) \cong B(\lambda ),$
\end{center}
where $B(\lambda )$ denotes the crystal basis of the integrable highest 
weight $U_q (\mathfrak{g})$-module $V(\lambda )$ with highest weight $\lambda $.
\label{1.0.2}
\end{thm}
Consequently, we obtain an embedding:
\begin{center}
$B(\lambda )\hookrightarrow \widetilde{B}(\Tilde{\lambda })$,
\end{center}
where $\widetilde{B}(\Tilde{\lambda })$ denotes the crystal basis of the 
integrable highest weight $U_q (\Tilde{\mathfrak{g}})$-module $\widetilde{V}(\Tilde{\lambda })$
with highest weight $\Tilde{\lambda }$. This embedding is not necessarily a morphism of crystals.
However, it is a \textit{quasi-embedding} of crystals in the sense of \S 4.1 below;
for the precise statement, see \S 4 and \S 6.

By using these results, we obtain the following decomposition rules.
\begin{thm}$(\mathrm{Decomposition\ Rule\ for\ tensor\ products}).$\ 
Let $\lambda , \mu \in P^+$, then we have
$$\mathbb{B}(\lambda )\otimes \mathbb{B}(\mu )\ \cong \ 
\bigsqcup _
{\begin{subarray}{c}
\pi \in \mathbb{B}(\mu )\\ 
\Tilde{\pi }\ \!\!:\ \!\! \Tilde{\lambda }
\text{-}\mathrm{dominant}
\end{subarray}}
\mathbb{B}\bigl( \lambda +\pi (1) \bigl).$$
Here, $\Tilde{\pi } \in \widetilde{\mathbb{B}}(\Tilde{\mu })$ denotes the 
image of $\pi \in \mathbb{B}(\mu )$ under the embedding 
$\mathbb{B}(\mu ) \hookrightarrow \widetilde{\mathbb{B}}(\Tilde{\mu })$,
and it is said to be $\Tilde{\lambda }$-dominant
if $\Tilde{\pi }(t)+\Tilde{\lambda }$ belongs to the dominant Weyl chamber 
of $\Tilde{\mathfrak{g}}$ for all $t\in [0,1]$.
\end{thm}
\begin{thm}$(\mathrm{Branching\ Rule\ for\ restriction\ to\ Levi\ subalgebras}).$\ 
Let $\lambda \in P^+$. For a subset $S \subset I$, 
we denote by $\mathfrak{g}_S$ the corresponding Levi subalgebra 
of $\mathfrak{g}$, by $\Tilde{\mathfrak{g}}_{\widetilde{S}}$ the one for $\Tilde{\mathfrak{g}}$
(see \S 7.2), and by $\mathbb{B}_S (\lambda )$ the set of GLS paths of shape $\lambda $ for $\mathfrak{g}_S$. Then,
$$\mathbb{B}(\lambda )\ \cong \ \bigsqcup _
{\begin{subarray}{c}\pi \in \mathbb{B}(\lambda) \\ 
\Tilde{\pi }\ \!\! :\ \!\! \Tilde{\mathfrak{g}}_{\widetilde{S}}
\text{-}\mathrm{dominant}
\end{subarray}} \mathbb{B}_S \bigl( \pi (1) \bigl) \ \ 
\mathrm{as}\ \mathfrak{g}_S \mathrm{\text{-}crystals}.$$
Here, $\Tilde{\pi } \in \widetilde{\mathbb{B}}(\Tilde{\lambda })$ 
denotes the image of $\pi \in \mathbb{B}(\lambda )$ under the 
embedding $\mathbb{B}(\lambda ) \hookrightarrow \widetilde{\mathbb{B}}(\Tilde{\lambda })$,
and it is said to be $\Tilde{\mathfrak{g}}_{\widetilde{S}}$-dominant
if $\Tilde{\pi }(t)$ belongs to the dominant Weyl chamber of $\Tilde{\mathfrak{g}}_{\widetilde{S}}$
for all $t\in [0,1]$.
\end{thm}
Note that these theorems extend the corresponding (well-known) results for Kac--Moody algebras
(see [\textbf{Kas1-2,4, Li1-2}]).

This paper is organized as follows.
In Section 2, we recall some elementary facts about generalized Kac--Moody algebras, and 
introduce the monoid $\mathcal{W}$.
In Section 3, we review Joseph--Lamprou's path model.
In Section 4, we introduce our main tool, which we call an embedding of path models, 
and then give a new proof of the isomorphism theorem for path crystals.
In Section 5, we give a characterization of standard paths. 
In Section 6, we prove that the crystal of GLS paths of shape $\lambda $
is isomorphic to the crystal basis of the integrable highest weight module with highest weight $\lambda $.
In Section 7, we obtain some decomposition rules for path crystals for generalized Kac--Moody algebras.
In the Appendix, we give the (postponed) proofs of results on the 
monoid $\mathcal{W}$ stated in Section 2.
\begin{flushleft}
\textbf{Acknowledgments}
\end{flushleft}
This paper is a part of the author's master thesis [\textbf{I}], which was written under the supervision of 
Professor Satoshi Naito. The author would like to express his sincere thanks 
to Professor Satoshi Naito and Professor Daisuke Sagaki for suggesting 
the problems taken up in this paper and for many helpful discussions. 

\section{Preliminaries}

\subsection{Generalized Kac--Moody algebras}
In this subsection, we recall some fundamental facts about 
generalized Kac--Moody algebras. For more details, 
we refer the reader to [\textbf{Bo1-2, JKK, JKKS, JL, Kac, Kan}].

Let $I$ be a countable index set. We call 
$A=(a_{ij})_{i,j\in I}$ a \textit{Borcherds--Cartan\ matrix}
if the following three conditions are satisfied:
(1)\ $a_{ii}=2\ \mathrm{or} \ 
a_{ii}\in \mathbb{Z}_{\le 0}\ 
\mathrm{for\ each}\ i\in I;$
(2)\ $a_{ij}\in \mathbb{Z}_{\le 0}\ 
\mathrm{for\ all}\ i,j\in I,\ \mathrm{with}\ i\neq j;$
(3)\ $a_{ij}=0\ \mathrm{if\ and\ only\ if}\ a_{ji}=0\ 
\mathrm{for\ each}\ i,j\in I,\ \mathrm{with}\ i\neq j.$

An index $i\in I$ is said to be \textit{real} if $a_{ii}=2$, and 
\textit{imaginary} if $a_{ii}\le 0$. Denote by 
$I^{re}:=\{ i\in I\ |\ a_{ii}=2 \}$ the set of real indices, and by 
$I^{im}:=\{ i\in I\ |\ a_{ii}\le 0 \} =I\setminus I^{re}$ the set of 
imaginary indices. A Borcherds--Cartan matrix $A$ is said to be 
\textit{symmetrizable} if there exists a diagonal matrix 
$D=\mathrm{diag}(d_i ) _{i\in I}$, with $d_i \in \mathbb{Z}_{>0}$, 
such that $DA$ is symmetric. Also, if $a_{ii}\in 2\mathbb{Z}$ for all $i\in I$, then $A$ is said to be \textit{even}. 
In \S 6 and \S 7, we assume that the Borcherds--Cartan matrix is symmetrizable and even. Until then, we take an arbitrary Borcherds--Cartan matrix.

For a given Borcherds--Cartan matrix $A=(a_{ij})_{i,j\in I}$, 
a \textit{Borcherds--Cartan datum} is a quintuple 
$\bigl( A,\varPi :=\{ \alpha _i \} _{i\in I}, \varPi ^{\vee}:=\{ \alpha _i ^{\vee} \} _{i\in I},
P, P^{\vee}\bigl) ,$ where $\varPi $ and $\varPi ^{\vee}$ are the sets of
\textit{simple roots} and \textit{simple coroots}, respectively, $P^{\vee}$ 
is a \textit{coweight lattice}, and $P:=\mathrm{Hom}_{\mathbb{Z}}(P^{\vee},\mathbb{Z})$ 
is a \textit{weight lattice}. We set $\mathfrak{h}:=P^{\vee}\otimes _{\mathbb{Z}}\mathbb{C}$,
and call it the \textit{Cartan subalgebra}. Let
$\mathfrak{h}^* :=\mathrm{Hom}_{\mathbb{C}}(\mathfrak{h},\mathbb{C})$ 
denote the full dual space of $\mathfrak{h}$.
In this paper, we assume that $\varPi ^{\vee} \subset \mathfrak{h}$
and $\varPi \subset \mathfrak{h}^*$ are both linearly independent over $\mathbb{C}$.
Let $P^+ :=\{ \lambda \in P \ |\ \alpha _i ^{\vee}(\lambda )
\ge 0 \ \mathrm{for\ all}\ i\in I \}$ be the set of dominant integral weights, and 
$Q:=\bigoplus _{i\in I}\mathbb{Z}\alpha _i $ the \textit{root lattice}; 
we set $Q^+ :=\sum _{i\in I}\mathbb{Z}_{\ge 0}\alpha _i $.
Also, we write $\mu \succ \nu $ if $\mu -\nu \in Q^+$.
Let $\mathfrak{g}$ be the generalized Kac--Moody algebra associated with 
a Borcherds--Cartan datum $(A,\varPi ,\varPi ^{\vee},P,P^{\vee} )$. 
We have the root space decomposition
$\mathfrak{g}=\bigoplus _{\alpha \in \mathfrak{h}^*}
\mathfrak{g}_{\alpha }$, where $\mathfrak{g}_{\alpha }
:=\{ x\in \mathfrak{g}\ |\ \mathrm{ad}(h)(x)=h(\alpha )x \ 
\mathrm{for\ all}\ h\in \mathfrak{h} \},$ and $\mathfrak{h}=\mathfrak{g}_0.$
Denote by $\varDelta :=\{ \alpha \in \mathfrak{h}^* \mid 
\mathfrak{g}_{\alpha }\neq \{ 0 \} ,\ \alpha \neq 0 \}$ 
the root system, and by $\varDelta ^+$ the set of positive roots.
Note that $\varDelta \subset Q \subset P$ and 
$\varDelta =\varDelta ^+ \sqcup (- \varDelta ^+)$.
Let $r_i \in \mathrm{GL}(\mathfrak{h}^*)$, $i\in I^{re}$,  
denote the simple reflections defined by 
$r_i (\mu) := \mu - \alpha _i ^{\vee}(\mu ) \alpha _i $ for $\mu \in \mathfrak{h}^*$.
The subgroup $\mathcal{W}_{re} \subset \mathrm{GL}(\mathfrak{h}^*)$ 
generated by the simple reflections $r_i$, $i\in I^{re}$, is a Coxeter group.
Set $\varDelta _{im}:=\mathcal{W}_{re}\varPi _{im}$,
where $\varPi _{im}:=\{ \alpha _i \} _{i\in I^{im}}$;
note that $\varDelta _{im} \subset \varDelta ^+$.
Also, we set $\varDelta _{re}=\mathcal{W}_{re}\varPi _{re}$,
where $\varPi _{re}:=\{ \alpha _i \} _{i\in I^{re}}$, and 
$\varDelta _{re}^+ :=\varDelta _{re} \cap \varDelta ^+$; note 
that $\varDelta _{re}$ is the set of real roots.
As in the case of ordinary Kac--Moody algebras, 
it is easily checked that the coroot 
$\beta ^{\vee}:=w\alpha _i ^{\vee}$ of
$\beta =w\alpha _i \in \varDelta _{re}^+ \sqcup \varDelta _{im}$
is well-defined (see [\textbf{JL}, \S 2.1.9]).

\subsection{The monoid $\mathcal{W}$ and its properties}
In this subsection, we associate a certain monoid $\mathcal{W}$ 
to a given Borcherds--Cartan datum. For an ordinary Cartan datum, 
namely, a Borcherds--Cartan datum with $I^{im} = \emptyset$, 
the associated monoid $\mathcal{W}$ is in fact a Weyl group of this Cartan datum. 
So, this monoid can be thought of as a generalization of a Weyl group.
Also, we obtain some combinatorial properties of this monoid.
Since several results stated in this subsection are not used explicitly in this paper, 
we postpone giving the proofs of these results to the Appendix.

Let $\bigl( A=(a_{ij}) _{i,j\in I}, \varPi =\{ \alpha _i \} _{i\in I}, \varPi ^{\vee}=\{ \alpha _i ^{\vee} \} _{i\in I},
P, P^{\vee}\bigl)$ be a Borcherds--Cartan datum as in the previous section.  
If we define $r_i \in \mathrm{GL}(\mathfrak{h}^*)$, $i\in I^{im}$, by
$r_i (\mu ):=\mu -\alpha _i ^{\vee}(\mu )\alpha _i $ for $\mu \in \mathfrak{h}^*$,
then the inverse of $r_i $ is as follows: 
$$r_i ^{-1}(\mu )=\mu +\frac{1}{1-a_{ii}}\alpha _i ^{\vee}(\mu )\alpha _i\ 
\mathrm{for}\ \mu \in \mathfrak{h}^*.$$
Note that this $r_i$ has an infinite order in $\mathrm{GL}(\mathfrak{h}^*)$.
Still, we can verify that the element $r_{\beta }:=wr_i w^{-1}\in \mathrm{GL}(\mathfrak{h}^*)$ 
for $\beta =w\alpha _i \in \varDelta _{re}^+ \sqcup \varDelta _{im}$ 
and $w\in \mathcal{W}_{re}$, is well-defined.
By convention, we call $r_{\beta }$ the reflection with respect to the root $\beta $.
\begin{define}\label{2.2.1}
Let $\mathcal{W}$ denote the monoid generated
by the elements $\Tilde{r}_i ,\ i\in I$, 
subject to the following relations:
\begin{itemize}
 \item[(1)]
$\Tilde{r}_i ^2 =1$ for all $i\in I^{re}$;
 \item[(2)]
if $i,j \in I^{re},\ i\neq j,$ and the order of 
$r_i r_j \in \mathrm{GL}(\mathfrak{h}^*)$ is 
$m\in \{ 2,3,4,6 \}$, then we have 
\begin{center}
$(\Tilde{r}_i \Tilde{r}_j )^m 
=(\Tilde{r}_j \Tilde{r}_i )^m =1$;
\end{center}
 \item[(3)]
let $i\in I^{im}$. Then for all $j\in I\setminus \{ i\}$ such that 
$a_{ij}=0$, we have $\Tilde{r}_i \Tilde{r}_j =\Tilde{r}_j \Tilde{r}_i$.
\end{itemize}
\end{define}
Each element $w \in \mathcal{W}$ can be written as 
a product $w= \Tilde{r}_{i_k} \cdots \Tilde{r}_{i_2} \Tilde{r}_{i_1}$ 
of generators $\Tilde{r}_i$, $i \in I$. If the number $k$ is minimal among 
all the expressions for $w$ of the form above, then $k$ is called the 
\textit{length} of $w$ and the expression 
$\Tilde{r}_{i_k} \cdots \Tilde{r}_{i_2} \Tilde{r}_{i_1}$
is called a \textit{reduced\ expression}. In this case, 
we write $\ell _{\mathcal{W}}(w) = \ell (w) :=k$.
Since the $r_i \in \mathrm{GL}(\mathfrak{h}^*)$, $i\in I$, satisfy the conditions 
$(1), (2)$ and $(3)$ of Definition \ref{2.2.1},  
we have the following (well-defined) homomorphism of monoids:
\begin{center}
$\mathcal{W}\longrightarrow \mathrm{GL}(\mathfrak{h}^*),
\ \Tilde{r}_i \longmapsto r_i \ (i\in I).$
\end{center}
By abuse of notation, we also write $r_i $ for $\Tilde{r}_i$ in $\mathcal{W}$.
\begin{conj}\label{2.2.2}
The homomorphism 
$\mathcal{W}\longrightarrow \mathrm{GL}(\mathfrak{h}^*)$ defined above is injective.
\end{conj}

In order to study properties of the monoid $\mathcal{W}$,
we introduce a Coxeter group $\mathfrak{W}$.
\begin{define}\label{2.2.3} 
We set $\Tilde{I}:=\{ (i,1) \} _{i\in I^{re}} \sqcup \{ (i,m) \} _{i\in I^{im}, m\in \mathbb{Z}_{\ge 1}}
\subset I \times \mathbb{Z}$. Let $\mathfrak{W}$ denote the Coxeter group associated with 
the Coxeter matrix $X=\bigl( \mathrm{x}_{(i,m), (j,n)} \bigl) 
_{(i,m), (j,n) \in \Tilde{I}}$ given by:
\begin{itemize}
\item[(1)]
$\mathrm{x}_{(i,m), (i,m)} =1$ for all $(i,m) \in \Tilde{I}$;
\item[(2)]
if $i,j\in I^{re}$, then $\mathrm{x}_{(i,1), (j,1)}$
is the order of $r_i r_j$ in $\mathrm{GL}(\mathfrak{h}^*)$;
\item[(3)]
let $i\in I^{im}$ and $m\in \mathbb{Z}_{\ge 1}$.
Then for all $(j,n) \in \Tilde{I}\setminus \{ (i,m) \}$, we have 
$$\mathrm{x}_{(i,m), (j,n)}=\mathrm{x}_{(j,n), (i,m)}=
\begin{cases}
\ 2 & \mathrm{if}\ a_{ij}=0, \\
\ \infty & \mathrm{otherwise}.
\end{cases}$$
\end{itemize}
\end{define}
Let $s_{(i,m)}$, $(i,m) \in \Tilde{I}$, denote the simple reflections of $\mathfrak{W}$, 
and let $\ell _{\mathfrak{W}} : \mathfrak{W} \longrightarrow \mathbb{Z}_{\ge 1}$ denote the length function 
of $\mathfrak{W}$ (see [\textbf{Hu}, \S 5.2]).

For notational simplicity, we use ``multi-index" notation.
We write $(\mathbf{i,m}) =((i_s , m_s))_{s=1} ^k
=((i_k , m_k),\ldots , (i_2, m_2), (i_1, m_1))$ if 
$\mathbf{i}=(i_k ,\ldots ,i_2 ,i_1) \in I^k$ and 
$\mathbf{m}=(m_k ,\ldots , m_2 , m_1) \in \mathbb{Z}^k$
for $k\ge 0$, and call it an \textit{ordered index} if 
$m_s =1$ for $i_s \in I^{re}$, and if $m_{x_q} =q$ for all 
$q=1,2,\ldots , p$, where $\{ x_1 , x_2 , \ldots , x_p \} 
=\{ 1\le  x \le k \mid i_x = i \}$, with $1 \le x_1 < x_2 < \cdots < x_p \le k$, for $i\in I^{im}$.
We set $\mathcal{I}:=\bigcup _{k=1} ^{\infty} I^k$, 
$\widetilde{\mathcal{I}}:=\bigcup _{k=1}^{\infty}
\Tilde{I}^k$, and denote by $\widetilde{\mathcal{I}}_{\mathrm{ord}}$ the set of 
all ordered indices. Note that if $\mathbf{(i,m)} \in \widetilde{\mathcal{I}}_{\mathrm{ord}}$,
then $\mathbf{m}$ is determined uniquely by $\mathbf{i}$.
Therefore, we have a bijection 
$\mathcal{I}\rightarrow \widetilde{\mathcal{I}}_{\mathrm{ord}} ,\ 
\mathbf{i}\mapsto \mathbf{(i,m)}.$
\begin{ex}\label{2.2.4}
If $I = \{ 1, \bm{2}, \bm{3} \}$, $I^{re} =\{ 1 \}$ and $I^{im} = \{ \bm{2}, \bm{3} \}$, 
then $\Tilde{I} = \{ (1,1) \} \sqcup \{ ( \bm{2}, m) ,$
$( \bm{3}, m) \} _{m\in \mathbb{Z}_{\ge 1}}$.
An ordered index $\mathbf{(i,m)} \in \widetilde{\mathcal{I}}_{\mathrm{ord}}$
corresponding to 
$\mathbf{i}=(\bm{3},1,\bm{3},\bm{3},1,\bm{2},\bm{3},1,\bm{3},\bm{2},1,\bm{2})\in \mathcal{I}$
is as follows:
\begin{center}
$\mathbf{(i,m)}=\bigl( (\bm{3},5),(1,1),(\bm{3},4),(\bm{3},3),(1,1),(\bm{2},3),(\bm{3},2),(1,1),
(\bm{3},1),(\bm{2},2),(1,1),(\bm{2},1) \bigl)$.
\end{center}
\end{ex}

Now, we define a subset $\mathfrak{V} \subset 
\mathfrak{W}$ by $\mathfrak{V}:=\{ S_{\mathbf{(i,m)}} \in \mathfrak{W} \mid 
\mathbf{(i,m)} \in \widetilde{\mathcal{I}}_{\mathrm{ord}} \} ,$
where $S_{\mathbf{(i,m)}}:=
s_{(i_k , m_k)} \cdots s_{(i_2 , m_2)} 
s_{(i_1, m_1)}$ with $\mathbf{i} =(i_k ,\ldots ,i_2 ,i_1)$ 
and $\mathbf{m}=(m_k ,\ldots , m_2 , m_1)$; 
note that $\mathfrak{V}$ is not a subgroup of $\mathfrak{W}$
since it is not closed under multiplication.
Also, we write $R_{\mathbf{i}}:=r_{i_k }\cdots r_{i_2} r_{i_1 } \in \mathcal{W}$ 
with $\mathbf{i}=(i_k ,\ldots ,i_2 ,i_1) \in \mathcal{I}$.
\begin{lem}\label{2.2.5}
The map (not necessarily a homomorphism of monoids) given below is bijective.
$$\sigma : \mathcal{W}\longrightarrow 
\mathfrak{V},\ R_{\mathbf{i}}\longmapsto 
S_{\mathbf{(i,m)}} ,$$
where the $\mathbf{m}$ for which $\mathbf{(i,m)}\in \widetilde{\mathcal{I}}_{\mathrm{ord}}$
is determined uniquely by $\mathbf{i}$.
\end{lem}
Observe that $\ell _{\mathcal{W}} (w)=\ell _{\mathfrak{W}} \bigl( \sigma (w) \bigl)$ for each $w\in \mathcal{W}$, 
and that $\mathcal{W}_{re}$ (see \S 2.1) and 
$\langle \Tilde{r}_i \ |\ i\in I^{re} \rangle _{\mathrm{monoid}} \subset \mathcal{W}$
are isomorphic as groups, where $\langle \Tilde{r}_i \ |\ i\in I^{re} \rangle _{\mathrm{monoid}}$ denotes 
the submonoid of $\mathcal{W}$ generated by $\Tilde{r}_i $, $i\in I^{re}$; 
note that this submonoid $\langle \Tilde{r}_i \ |\ i\in I^{re} \rangle _{\mathrm{monoid}}$
is in fact a group since each $\Tilde{r}_i $, $i\in I^{re}$, is an involution by Definition \ref{2.2.1} (1). 
Hence we may (and do) regard the group
$\mathcal{W}_{re}$ as a submonoid of $\mathcal{W}$.
Also, it is clear that the correspondence, 
$\mathcal{W} \ni w \Tilde{r}_i w^{-1} \mapsto wr_i w^{-1} \in \mathrm{GL}(\mathfrak{h}^*)$
for $w\in \mathcal{W}_{re}$ and $i\in I$, is injective.
For this reason, we write $r_{w(\alpha _i)}$ for $wr_i w^{-1}$ in $\mathcal{W}$.
\begin{define}$\label{2.2.6}
\mathrm{(cf.\ [\mathbf{BB,\ Hu}]).}$\ 
For $w\in \mathcal{W}$ and 
$\beta \in \varDelta _{re} ^+ \sqcup \varDelta _{im}$, we write
$w \rightarrow r_{\beta} w$ or $w\xrightarrow{\beta }r_{\beta }w$
if $\ell (r_{\beta } w)>\ell (w)$. Also, we define a partial order
$\le $ on $\mathcal{W}$ as follows: $w \le w' \ \mathrm{in}\ \mathcal{W}$ if there exist 
$w_0 ,w_1 ,\ldots ,w_l \in \mathcal{W}$
such that $w=w_0\rightarrow w_1 \rightarrow \cdots \rightarrow w_l =w' .$
We call this partial order on $\mathcal{W}$ the Bruhat order.
\end{define}
Note that $w \rightarrow r_{\beta }w$ for all 
$w\in \mathcal{W}$ and $\beta \in \varDelta _{im}$.

As a corollary of Lemma \ref{2.2.5}, we can show that 
the \textit{Word Property} (cf. [\textbf{BB}, Theorem 3.3.1]) holds for $\mathcal{W}$.
Also, the following \textit{Exchange Property} (cf. [\textbf{BB,\ Hu}]) holds for $\mathcal{W}$.
\begin{cor}$\label{2.2.7}
\mathrm{(Strong\ Exchange\ Property\ for\ real\ roots;\ cf.\ [\mathbf{BB,\ Hu}]).}$\ 
Let $w=r_{i_k}\cdots r_{i_1}$
$\in \mathcal{W}$.
If there exists $\beta \in \varDelta _{re} ^+$ such that
$\ell (r_{\beta }w)<\ell (w)$, then there exists some $i_s \in I^{re}$, 
$1\le s \le k$, such that $r_{\beta }w=r_{i_k}\cdots \widehat{r_{i_s}}
\cdots r_{i_1}$, where $\widehat{r_{i_s}}$ indicates that the $r_{i_s}$ is omitted 
in the expression. Moreover, if $r_{i_k}\cdots r_{i_1}$ is a reduced expression, 
then the $i_s \in I^{re}$ above is uniquely determined.
\end{cor}
\begin{cor}\label{2.2.8}
$\mathrm{(Strong\ Exchange\ Property\ for\ imaginary\ roots;\ cf.\ [\mathbf{BB,\ Hu}]).}$ 
Let $\beta =w \alpha _i$
$\in \varDelta _{im},\ i\in I^{im},\ w\in \mathcal{W}_{re},$ and let
$v, v' \in \mathcal{W}$. If $v' \xrightarrow{\beta } v$,
then an expression of $v'$ is obtained from the expression $v=r_{i_k }\cdots r_{i_1}$
by omitting the leftmost $r_i$, i.e., if 
$\{ x_1,\ldots , x_p \} = \{ 1\le x \le k \mid i_x =i \}$ with $1\le x_1 <x_2 <\cdots <x_p \le k$,
then $v' =r_{i_k }\cdots \widehat{r_{i_{x_p}}} \cdots r_{i_{x_2}} \cdots r_{i_{x_1}} \cdots r_{i_1}$.
\end{cor}

By an argument similar to the one for [\textbf{JL}, Lemma 2.2.3], 
we can show that $\mathcal{W}$ has the following property.
\begin{lem}$(\mathrm{cf.\ [\mathbf{JL},\ Lemma\ 2.2.3}]).$\ 
Let $\lambda \in P^+$ and $i\in I^{im}$.
If there exists $1\neq w \in \mathcal{W}_{re}$
such that $r_i w \lambda \notin P^+$, then
there exists $j\in I^{re}$ such that 
$r_i w \lambda =r_j r_i w' \lambda $,
where $w' =r_j w$ with $\ell (w')=\ell (w)-1$,
and such that $r_i r_j =r_j r_i $ in $\mathcal{W}$.
Therefore, there exist $w_1 ,\ w_2 \in \mathcal{W}_{re}$
such that $w=w_1 w_2$ with $\ell (w)=\ell (w_1)+\ell (w_2)$,
and such that $r_i w_1 =w_1 r_i$ in $\mathcal{W}$, $r_i w_2 \lambda \in P^+$.
\label{2.2.9}
\end{lem}
\begin{define}\label{2.2.10}
$(\mathrm{cf.\ [\mathbf{JL},\ \S 2.2.6]).}$\ 
Let $w=w_k r_{i_k} \cdots w_1 r_{i_1} w_0 \in \mathcal{W}$ 
be a reduced expression, where 
$i_1 ,\ldots ,i_k \in I^{im}$ and $w_0 ,w_1 ,\ldots ,w_k \in \mathcal{W}_{re}$. 
We call this expression a dominant reduced expression if it satisfies 
\begin{center}
$r_{i_s} w_{s-1} r_{i_{s-1}} \cdots w_1 r_{i_1} w_0 (P^+) \subset P^+$
\end{center}
for all $s=1,2,\ldots ,k$. 
\end{define}
Note that every $w\in \mathcal{W}$ 
has at least one dominant reduced expression.
Indeed, if we choose an expression of the form above in such a way that
the sequence $\bigl( \ell (w_0 ), \ell (w_1), \ldots , \ell (w_k ) \bigl)$ 
is minimal in lexicographic order among all the reduced expressions 
of the form above, then it is a dominant reduced expression by Lemma \ref{2.2.9}.
\begin{lem}$(\mathrm{cf.\ [\mathbf{JL},\ Lemma\ 2.2.2]}).$\ 
For every $\lambda \in P^+$, we have
$\mathrm{Stab}_{\mathcal{W}}(\lambda ) 
:=\{ w\in \mathcal{W}\ |$
$w(\lambda )=\lambda \}
=\langle r_i \ | \ \alpha _i ^{\vee}(\lambda )=0 
\rangle _{\mathrm{monoid}}$, where the right-hand side is the monoid 
generated by those $r_i ,\ i\in I$, for which $\alpha _i ^{\vee}(\lambda )=0$.
\label{2.2.11}
\end{lem}

\section{Joseph--Lamprou's path model}
In this section, following [\textbf{JL}], we review Joseph--Lamprou's path model
for generalized Kac--Moody algebras.

\subsection{Root operators}
Let $\mathfrak{h}_{\mathbb{R}}$ denote a real form of $\mathfrak{h}$,
and $\mathfrak{h}_{\mathbb{R}} ^*$ be its full dual space.
Let $\mathbb{P}$ be the set of all piecewise-linear continuous maps
$\pi :[0,1]\longrightarrow \mathfrak{h}_{\mathbb{R}} ^*$
such that $\pi (0)=0$ and $\pi (1) \in P$, 
where we set $[0,1]:=\{ t\in \mathbb{R} \mid 0\le t\le1 \}$. 
Also, we set $H_i ^{\pi }(t):=\alpha _i ^{\vee}\bigl( \pi (t) \bigl)$ for 
$t\in [0,1]$, and then $m_i ^{\pi }:=\min \{ H_i ^{\pi }(t) \mid H_i ^{\pi }(t)\in \mathbb{Z},\ t\in [0,1]\}$.
Remark that we do not identify two paths $\pi _1 , \pi _2 \in \mathbb{P}$ even if 
there exist piecewise linear, nondecreasing, surjective, continuous maps 
$\psi _1 , \psi _2 :[0,1]\longrightarrow [0,1]$ such that $\pi _1 \circ \psi _1 =\pi _2 \circ \psi _2$.

Now we define the root operators 
$e_i ,f_i :\mathbb{P}\longrightarrow \mathbb{P} \sqcup \{ \mathbf{0} \}$ for $i\in I$.
First, we set 
$f_+ ^i (\pi ):= \max \{ t\in [0,1] \mid H_i ^{\pi }(t)=m_i ^{\pi} \}$. If $f_+ ^i (\pi )<1$, then we can define 
$f_- ^i (\pi ) := \min \{ t\in [f_+ ^i (\pi ),1] \mid H_i ^{\pi }(t)=m_i ^{\pi }+1 \}$. In this case, we set 
\begin{align*}
(f_i \pi )(t):=
\begin{cases}
\pi (t) & t\in [0,f_+ ^i (\pi )], \\
\pi \bigl ( f_+ ^i (\pi ) \bigl ) + r_i \bigl ( \pi (t) -
\pi \bigl ( f_+ ^i (\pi ) \bigl ) \bigl ) \\
\hspace{1in}
=\pi (t) - \bigl( H_i ^{\pi} (t) - m_i ^{\pi} \bigl) \alpha _i
& t\in [f_+ ^i (\pi ), f_- ^i (\pi )], \\
\pi (t)-\alpha _i & t\in [f_- ^i (\pi ),1].
\end{cases}
\end{align*}
Otherwise (i.e., if $f_+ ^i (\pi )=1$), we set $f_i \pi :=\mathbf{0}$.

Next, we define the operator $e_i $ for $i\in I^{re}$.
Set $e _+ ^i (\pi ):= \min \{ t\in [0,1] \mid H_i ^{\pi }(t) 
=m_i ^{\pi } \}$. If $e_+ ^i (\pi )>0$, then we can define
$e_- ^i (\pi ):= \max \{ t\in [0,e_+ ^i (\pi )] \mid 
H_i ^{\pi }(t)= m_i ^{\pi } +1 \}$. In this case, we set
\begin{align*}
(e_i \pi )(t):=
\begin{cases}
\pi (t) & t\in [0,e_- ^i (\pi )], \\
\pi \bigl ( e_- ^i (\pi ) \bigl ) + r_i \bigl ( \pi (t) -
\pi \bigl ( e_- ^i (\pi ) \bigl ) \bigl ) \\
\hspace{1in}
=\pi (t) - \bigl( H_i ^{\pi}(t) - m_i ^{\pi} -1 \bigl) \alpha _i
& t\in [e_- ^i (\pi ), e_+ ^i (\pi )], \\
\pi (t)+\alpha _i & t\in [e_+ ^i (\pi ),1].
\end{cases}
\end{align*}
Otherwise (i.e., if $e_+ ^i (\pi )=0$), we set $e_i \pi :=\mathbf{0}$.

Finally, we define the operator $e_i $ for $i\in I^{im}$.
Set $e_- ^i (\pi ):=f_+ ^i (\pi )$. If $e_- ^i (\pi )<1$
and there exists $t \in [e_- ^i (\pi ),1]$ such that 
$H_i ^{\pi }(t)\ge m_i ^{\pi }+1-a_{ii}$, then we can define
$e_+ ^i (\pi ):=\min \{ t\in [e_- ^i (\pi ),1] \mid 
H_i ^{\pi }(t)=m_i ^{\pi }+1 -a_{ii} \}$. We set $e_i \pi :=\mathbf{0}$ if $e_- ^i (\pi )=1$, or $e_- ^i (\pi )<1$ and 
$H_i ^{\pi }(t) < m_i ^{\pi } + 1 - a_{ii} $ for all 
$t\in [e_- ^i (\pi ),1]$, or $e_- ^i (\pi )<1$ and   
there exists $t\in [e_- ^i (\pi ),1]$ such that 
$H_i ^{\pi }(t)\ge m_i ^{\pi }+1-a_{ii}$ and 
$H_i ^{\pi }(s) \le m_i ^{\pi }-a_{ii}$ for some 
$s\in [e_+ ^i (\pi ),1]$. Otherwise, we set
\begin{align*}
(e_i \pi )(t):=
\begin{cases}
\pi (t) & t\in [0,e_- ^i (\pi )], \\
\pi \bigl ( e_- ^i (\pi ) \bigl ) + r_i ^{-1} \bigl ( \pi (t) -
\pi \bigl ( e_- ^i (\pi ) \bigl ) \bigl ) \\
\hspace{1in} 
= \pi (t) + \frac{1}{1-a_{ii}}\bigl( H_i ^{\pi}(t) - m_i ^{\pi} \bigl) \alpha _i
& t\in [e_- ^i (\pi ), e_+ ^i (\pi )], \\
\pi (t)+\alpha _i & t\in [e_+ ^i (\pi ),1].
\end{cases}
\end{align*}
\begin{rem}\label{3.1.1}
From the definition of the root operator $f_i $, $i\in I$, 
there exists a continuous function 
$\psi : [0,1] \longrightarrow [0,1]$ which depends only on the 
function $H_i ^{\pi} (t)$ such that 
$\psi(t)=0$ for $t\in [0, f_+ ^i (\pi)]$, 
$\psi(t)=H_i ^{\pi}(t)-m_i ^{\pi} \in [0,1]$ for $t\in [f_+ ^i (\pi), f_- ^i (\pi)]$, 
$\psi(t)=1$ for $t\in [f_- ^i (\pi), 1]$, and 
such that $(f_i \pi )(t) =\pi (t) - \psi (t) \alpha _i $ for all $t \in [0,1]$
if $f_i \pi \neq \mathbf{0}$. Similar property holds for the 
root operator $e_i $, $i\in I$.
\end{rem}

\subsection{Generalized Lakshmibai--Seshadri paths}
Let $\lambda \in P^+$. We write 
$\mu \ge \nu $ for $\mu ,\nu \in \mathcal{W}\lambda
:= \{ w\lambda \in P \mid w\in \mathcal{W} \} $
if there exists a sequence of elements 
$\mu = \lambda _0 , \lambda _1 , \ldots , 
\lambda_k = \nu$ in $\mathcal{W} \lambda $
and positive roots $\beta _1 ,\ldots ,\beta _k 
\in \varDelta _{re}^+ \sqcup \varDelta _{im}$
such that $\lambda _{s-1} = r_{\beta _s} \lambda _s$
and $\beta _s ^{\vee}(\lambda _s )>0$ for $s=1,\ldots , k$.
This relation $\ge $ on $\mathcal{W} \lambda$ defines a partial order.
Note that $\lambda$ is minimum in $\mathcal{W}\lambda$ with respect to this partial order.
For $\mu , \nu \in \mathcal{W} \lambda $ and $\beta \in \varDelta _{re}^+ \sqcup \varDelta _{im}$, 
we write $\mu \xleftarrow{\beta } \nu$ if $\mu = r_{\beta } \nu$, $\beta ^{\vee}(\nu )>0 $,
and $\mu $ covers $\nu $ by this partial order;
in this case, $\beta $ is determined uniquely by $\mu $ and $\nu $, 
and we have $\mu - \nu =r_{\beta } \nu - \nu = - \beta ^{\vee}(\nu ) \beta \in - Q^+$.
Note that the direction of the arrow $\xleftarrow{\beta}$ defined above and 
that of the Joseph--Lamprou's one defined in [\textbf{JL}, \S 5.1.1] are opposite. 
From the definitions, it follows that the natural surjection 
$\mathcal{W}\longrightarrow \mathcal{W}\lambda $ 
preserves the partial orders on $\mathcal{W}$ and $\mathcal{W}\lambda $.
In [\textbf{JL}], Joseph--Lamprou generalized the notion of Lakshmibai--Seshadri paths
to the case of generalized Kac--Moody algebras as follows.
\begin{define}\label{3.2.1}
$\mathrm{([\mathbf{JL},\ \S 5.2.1]}).$\ 
For a rational number $a\in (0,1]$ and 
$\mu ,\nu \in \mathcal{W}\lambda $ with $\mu \ge \nu $,
an $a$-chain for $(\mu ,\nu )$ is a sequence 
$\mu =\nu _0 \xleftarrow{\beta _1 }\nu _1
\xleftarrow{\beta _2 } \cdots \xleftarrow{\beta _k }\nu _k =\nu $
of elements in $\mathcal{W}\lambda $ such that for each $s=1,2,\ldots ,k,\ 
a\beta _s ^{\vee}(\nu _s)\in \mathbb{Z}_{>0}$
if $\beta _s \in \varDelta _{re} ^+,$ and 
$a\beta _s ^{\vee}(\nu _s)=1$ if $\beta _s \in \varDelta _{im}$.
\end{define}
Note that if there exists an $a$-chain for $(\mu ,\nu )$, then we have $a(\mu -\nu )\in -Q^+$.
\begin{define}\label{3.2.2}
$\mathrm{([\mathbf{JL},\ \S 5.2.2]).}$\ 
Let $\bm{\lambda }=(\lambda _1 >\lambda _2 >\cdots >\lambda _k)$
be a sequence of elements in $\mathcal{W}\lambda $, and
$\bm{a}=(0=a_0 <a_1 <\cdots <a_k =1)$ a sequence of rational numbers.
Then, the pair $\pi =(\bm{\lambda }; \bm{a})$ is called a generalized 
Lakshmibai--Seshadri path (GLS path for short) of shape $\lambda $
if it satisfies the following conditions (called the chain condition):
$(\mathrm{i})$ there exists an $a_s$-chain for 
$(\lambda _s ,\lambda _{s+1})$ for each $s=1,2,\ldots ,k-1$;
$(\mathrm{ii})$ there exists a $1$-chain for $(\lambda _k ,\lambda ).$
\end{define}
In this paper, we regard the pair $\pi =(\bm{\lambda };\bm{a})$ 
as a path belonging to $\mathbb{P}$ by
$\pi (t):=\sum _{l=1}^{s-1}(a_l -a_{l-1})\lambda _l +(t-a_{s-1})\lambda _s$
for $a_{s-1}\le t\le a_s$ and $s=1,2,\ldots ,k.$
Note that $\pi (1)\in \lambda -Q^+$.
We denote by $\mathbb{B}(\lambda )$ the set of all GLS paths of shape $\lambda $.
\begin{rem}\label{3.2.3}
(1) Joseph--Lamprou's GLS paths for an ordinary Cartan datum (namely, $I^{im} = \emptyset$) 
are in fact identical with Littelmann's LS paths. 
(2) Note that the GLS path $\pi $ is an integral and monotone path (see $[\mathbf{JL},\ \S 5.3.7$--$\S 5.3.9]$), 
i.e., any local minimums of $H_i ^{\pi}(t)$ are integers, and $H_i ^{\pi }(t)$ is 
monotone in the intervals $[f_+ ^i (\pi ),f_- ^i (\pi )]$ and $[e_- ^i (\pi ),e_+ ^i (\pi )]$ 
for all $i\in I$. Hence, for a GLS path $\pi \in \mathbb{B}(\lambda )$ for $\lambda \in P^+$,
we see that $f_i \pi \neq \mathbf{0}$ in $\mathbb{B}(\lambda )$ if and only if 
$H_i ^{\pi }(1)\ge m_i ^{\pi }+1$; $e_i \pi \neq \mathbf{0}$ in $\mathbb{B}(\lambda )$ for $i\in I^{re}$
if and only if $m_i ^{\pi }\le -1$; and $e_i \pi \neq \mathbf{0}$ in $\mathbb{P}$ for $i\in I^{im}$
if and only if $H_i ^{\pi }(1)\ge m_i ^{\pi }+1-a_{ii}$.
\end{rem}

\subsection{Crystal structure on $\mathbb{B}(\lambda )$}
In this subsection, we define a crystal structure on $\mathbb{B}(\lambda )$. 
Let $\pi =(\bm{\lambda }; \bm{a}) =(\lambda _1 ,\ldots ,\lambda _k ;
\ \! 0= a_0 ,a_1 ,\ldots ,a_k =1) \in \mathbb{B}(\lambda )$ be a GLS path of shape $\lambda $.
In this paper, we allow repetitions in $\bm{\lambda }$ and $\bm{a}$ for convenience; namely, we assume that
$\lambda _1 \ge \cdots \ge \lambda _k$ and $0=a_0 \le a_1 \le \cdots \le a_k =1$.

We set $f_{\pm } ^i =f_{\pm } ^i (\pi )$ for $i\in I$. 
Then it follows that $f_+ ^i =a_p$ for some $0\le p\le k$.
Also, if $f_i \pi \in \mathbb{P}$, then $a_{q-1} \le f_- ^i \le a_q$
for some $p+1 \le q \le k$.
\begin{prop}$\mathrm{([\mathbf{JL},\ Proposition\ 6.1.2]).}$\ 
Fix $i\in I^{re}$, and let $\pi \in \mathbb{B}(\lambda )$
be as above. If $f_i \pi \in \mathbb{P}$, then
\begin{center}
$f_i \pi =(\lambda _1 ,\ldots ,\lambda _p ,r_i \lambda _{p+1} ,
\ldots ,r_i \lambda _q ,\lambda _q ,\ldots ,\lambda _k;\ \! a_0 ,\ldots ,
a_{q-1}, f_- ^i , a_q ,\ldots ,a_k)$,
\end{center}
and, moreover, $f_i \pi \in \mathbb{B}(\lambda )$.
\label{3.3.1}
\end{prop}
\begin{prop}\label{3.3.2}
$\mathrm{([\mathbf{JL},\ Proposition\ 6.1.3]).}$\ 
Fix $i\in I^{im}$, and let $\pi \in \mathbb{B}(\lambda )$ be as above.
If $f_i \pi \in \mathbb{P}$, then $f^i _+ =0$ and 
\begin{center}
$f_i \pi =(r_i \lambda _1 ,\ldots ,r_i \lambda _q ,\lambda _q ,
\ldots ,\lambda _k;\ \! a_0 ,\ldots ,
a_{q-1}, f_- ^i , a_q ,\ldots ,a_k)$.
\end{center}
Moreover, in this case, $f_i \pi \in \mathbb{B}(\lambda )$
and $\alpha _i ^{\vee}(\lambda _1)=\alpha _i ^{\vee}(\lambda _2)=\cdots 
=\alpha _i ^{\vee}(\lambda _q) 
=\frac{1}{f_- ^i}$. Hence we have
$$f_i \pi =\Bigl( \lambda _1 -\frac{\alpha _i}{f_- ^i },
\ldots ,\lambda _q -\frac{\alpha _i}{f_- ^i },\lambda _q,
\ldots ,\lambda _k;\ \! a_0 ,\ldots ,
a_{q-1}, f_- ^i , a_q ,\ldots ,a_k \Bigl) .$$
\end{prop}
From Proposition \ref{3.3.2}, we see that the action of $f_i$, $i\in I^{im}$, on $\pi \in \mathbb{B}(\lambda )$
is determined completely by $f_- ^i$, and that $r_i $ (and $r_i ^{-1}$) commutes with the
reflections with respect to the positive roots appearing in $a$-chains of $\pi $ for $0<a \le f_- ^i$.

We set $e_{\pm } ^i =e_{\pm } ^i (\pi )$ for $i\in I$. Then it follows that 
$e_+ ^i =a_q$ for some $0\le q\le k$. Now, let $i\in I^{re}$.
If $e_i \pi \in \mathbb{P}$, then we have $a_{p-1} \le e_- ^i \le a_p$
for some $0 \le p \le q$.
\begin{prop}$\mathrm{([\mathbf{JL},\ Proposition\ 6.2.1]).}$\ 
Fix $i\in I^{re}$, and let $\pi \in \mathbb{B}(\lambda )$.
If $e_i \pi \in \mathbb{P}$, then
\begin{center}
$e_i \pi =(\lambda _1 ,\ldots ,\lambda _p ,r_i \lambda _p ,
\ldots , r_i \lambda _q ,\lambda _{q+1} ,\ldots ,\lambda _k ;\ \! 
a_0 ,\ldots , a_{p-1}, e_- ^i , a_p , \ldots , a_k)$,
\end{center}
and, moreover, $e_i \pi \in \mathbb{B}(\lambda )$.
\label{3.3.3}
\end{prop}
Next, we fix $i\in I^{im}$. Note that $e_i \pi \in \mathbb{P}$ 
does not necessarily imply $e_i \pi \in \mathbb{B}(\lambda )$. 
\begin{prop}\label{3.3.4}
Fix $i\in I^{im}$, and let $\pi \in \mathbb{B}(\lambda )$.
If $e_i \pi \in \mathbb{B}(\lambda )$, then $e^i _- =0$ and 
\begin{center}
$e_i \pi =(r_i ^{-1}\lambda _1 ,\ldots ,r_i ^{-1}\lambda _q ,
\lambda _{q+1} ,\ldots ,\lambda _k ;\ \! 
a_0 ,\ldots , a_{q-1} , e_+ ^i =a_q , \ldots , a_k).$
\end{center}
Moreover, in this case, $e_+ ^i \alpha _i ^{\vee}(\lambda _q )=1-a_{ii}$ and 
$\alpha _i ^{\vee}(\lambda _1)=\alpha _i ^{\vee}(\lambda _2)=\cdots 
=\alpha _i ^{\vee}(\lambda _q)=\frac{1-a_{ii}}{e_+ ^i}>0$.
Hence we have
$$e_i \pi =\Bigl( \lambda _1 +\frac{\alpha _i}{e_+ ^i },
\ldots ,\lambda _q +\frac{\alpha _i}{e_+ ^i },\lambda _{q+1},
\ldots ,\lambda _k ;\ \! a_0 ,\ldots ,
a_{q-1},e_+ ^i =a_q ,\ldots ,a_k \Bigl) .$$
\end{prop}
From Proposition \ref{3.3.4}, we see that the action of $e_i$, $i\in I^{im}$, on $\pi \in \mathbb{B}(\lambda )$
is determined completely by $e_+ ^i$, and that $r_i $ (and $r_i ^{-1}$) commutes with the reflections
with respect to the positive roots appearing in $a$-chains of $\pi $ for $0<a\le e_+ ^i =a_q$.

We only give the proof of Proposition \ref{3.3.4}; for the proofs of Propositions \ref{3.3.1}--\ref{3.3.3},
we refer the reader to [\textbf{JL},\ \S 6.1--\S 6.2]. For this purpose, we need the following refinement of 
[\textbf{JL}, Lemma 7.3.6]. For $\mu =\lambda -\sum _{i\in I}m_i \alpha _i 
\in \lambda -Q^+$, we set $\mathrm{depth}_i ^{\lambda }(\mu ):=m_i $.
\begin{lem}\label{3.3.5}
Let $\mu ,\nu \in \mathcal{W}\lambda $, and let $a\in (0,1]$ be a rational number, 
$\mu =\nu _0 \xleftarrow{\beta _1}\nu _1 \xleftarrow{\beta _2}\cdots 
\xleftarrow{\beta _k}\nu _k =\nu $ an $a$-chain for $(\mu ,\nu )$. 
If there exists $\beta _s =w\alpha _i $ with $1\le s \le k$, 
$i\in I^{im}$ and $w\in \mathcal{W}_{re}$ 
in this $a$-chain (that is, if $\mathrm{depth}_i ^0 (\mu -\nu )>0$)
such that $a\alpha _i ^{\vee}(\mu )\le 1-a_{ii}$,
then $\beta _s =\alpha _i $ and $a\alpha _i ^{\vee}(\mu )=1-a_{ii}$.
Therefore, $r_i ^{-1}\mu \in \mathcal{W}\lambda $, and there exists 
an $a$-chain for $( r_i ^{-1}\mu ,\nu )$.
\end{lem}
$\mathbf{Proof.}$\ 
Suppose the contrary, i.e., that $a\alpha _i ^{\vee}(\mu )<1-a_{ii}$.
In addition, suppose that $\beta _s =w\alpha _i \neq \alpha _i$.
If we write $\beta _s =w\alpha _i =\alpha _i +\gamma $
for some $\gamma \in \mathbb{Z}_{\ge 0}\varPi _{re}\setminus \{ 0\}$, 
then it follows that $\alpha _i ^{\vee}(\gamma )\le -1$, 
and hence $\alpha _i ^{\vee}(\beta _s )=a_{ii}+\alpha _i ^{\vee}(\gamma ) \le a_{ii}-1$.
Also, we have $1-a_{ii}>a\alpha _i ^{\vee}(\mu )=a\alpha _i ^{\vee}(\nu _s)
-a\beta _s ^{\vee}(\nu _s)\alpha _i ^{\vee}(\beta _s)-\sum _{l=1} ^{s-1}
a\beta _l ^{\vee}(\nu _l)\alpha _i ^{\vee}(\beta _l)\ge a\alpha _i ^{\vee}(\nu _s)
+(1-a_{ii})-\sum _{l=1} ^{s-1}a\beta _l ^{\vee}(\nu _l)\alpha _i ^{\vee}(\beta _l).$
Therefore, we obtain $0>a\alpha _i ^{\vee}(\nu _s)
-\sum _{l=1}^{s-1}a\beta _l ^{\vee}(\nu _l)\alpha _i ^{\vee}(\beta _l)\ge 0$,
which is a contradiction. Thus, we must have $\beta _s =\alpha _i $.
In this case, we have $1-a_{ii}>a\alpha _i ^{\vee}(\mu )=
a\alpha _i ^{\vee}(\nu _s)-a\alpha _i ^{\vee}(\nu _s)\alpha _i ^{\vee}(\alpha _i)-\sum _{l=1} ^{s-1}
a\beta _l ^{\vee}(\nu _l)\alpha _i ^{\vee}(\beta _l)=1-a_{ii}-\sum _{l=1} ^{s-1}a\beta _l ^{\vee}(\nu _l)
\alpha _i ^{\vee}(\beta _l),$ and hence $0>-\sum _{l=1} ^{s-1} a\beta _l ^{\vee}(\nu _l)
\alpha _i ^{\vee}(\beta _l) \ge 0$, which is again a contradiction.
Thus, we must have $a\alpha _i ^{\vee}(\mu )=1-a_{ii}$.
Consequently, it follows from [\textbf{JL}, Lemma 7.3.6] that $r_i ^{-1}\mu \in \mathcal{W}\lambda $
and that there exists an $a$-chain for $(r_i ^{-1}\mu ,\nu )$. \qed
\begin{flushleft}
\textbf{\textit{Proof of Proposition \ref{3.3.4}.}}
\end{flushleft}
Suppose the contrary, i.e., that $a_{q-1}< e_+ ^i <a_q$. Then, there exists an 
$e_+ ^i$-chain for $(r_i ^{-1} \lambda _q, \lambda _q)$.
In particular, $r_i ^{-1} \lambda _q >\lambda _q$.
Also, we have $r_i (r_i ^{-1}\lambda _q)=\lambda _q$, and 
$\alpha _i ^{\vee}(r_i ^{-1}\lambda _q)=\alpha _i ^{\vee}\bigl( \lambda _q +\frac{1}{1-a_{ii}}
\alpha _i ^{\vee}(\lambda _q)\alpha _i \bigl)=\frac{1}{1-a_{ii}}\alpha _i ^{\vee}(\lambda _q)>0.$
Therefore, we obtain $\lambda _q \xleftarrow{\alpha _i } r_i ^{-1} \lambda _q $,
which is a contradiction. Hence we conclude that $e_+ ^i =a_q $.
In this case, we have $a_q \alpha _i ^{\vee}(\lambda _q )\in \mathbb{Z}$
and $a_q \alpha _i ^{\vee}(\lambda _q )\le 1-a_{ii}$.
Indeed, if this is not the case, then $H_i ^{\pi }(t)=1-a_{ii}$ for some 
$t<e_+ ^i$ since $\alpha _i ^{\vee}(\lambda _1 )\ge \alpha _i ^{\vee}(\lambda _2)
\ge \cdots \ge \alpha _i ^{\vee}(\lambda _q )>0$, which is absurd.

Observe that we have the original $a_q (=e_+ ^i)$-chain 
$\lambda _q :=\nu _l \xleftarrow{\beta _{l-1}}\cdots 
\xleftarrow{\beta _1}\nu _1 =: \lambda _{q+1}$, 
and also a sequence $\lambda _q \xleftarrow{\alpha _i } r_i ^{-1}\lambda _q 
\leftarrow \cdots \leftarrow \lambda _{q+1}$,
which is obtained by combining $\lambda _q \xleftarrow{\alpha _i } r_i ^{-1}\lambda _q$ 
with the $e_+ ^i (=a_q)$-chain for $(r_i ^{-1}\lambda _q ,\lambda _{q+1})$
since we assume that $e_i \pi \in \mathbb{B}(\lambda )$.
From this observation, we see that $\mathrm{depth}_i ^0 (\lambda _q -\lambda _{q+1})>0$, and that
there exists some $\beta _s =w\alpha _i $, with $1\le s \le l-1$ and $w\in \mathcal{W}_{re}$.
Hence it follows from Lemma \ref{3.3.5} that 
$\beta _s =\alpha _i $ and $a_q \alpha _i ^{\vee}(\lambda _q )=1-a_{ii}$. \qed

\vspace{2mm}
Now, we define a crystal structure on $\mathbb{B}(\lambda )$. 
Let $\pi \in \mathbb{B}(\lambda )$. We set $\mathrm{wt}(\pi ):=\pi (1)\in \lambda -Q^+$.
For each $i\in I^{re}$, we set $\varepsilon _i (\pi ) := -m_i ^{\pi },\ 
\varphi _i (\pi ) := \alpha _i ^{\vee}\bigl( \mathrm{wt}(\pi )\bigl) -m_i ^{\pi }
=H_i ^{\pi }(1) -m_i ^{\pi }.$
Then, we have $\varepsilon _i (\pi )=\max \{ n\in \mathbb {Z}_{\ge 0} \ |\ e_i ^n \pi \in \mathbb{P} \}$, 
and $\varphi _i (\pi )=\max \{ n\in \mathbb {Z}_{\ge 0} \ |\ f_i ^n \pi \in \mathbb{P} \}$.
For each $i\in I^{im}$, we set $\varepsilon _i (\pi ):=0,\ 
\varphi _i (\pi ):= \alpha _i ^{\vee}\bigl( \mathrm{wt}(\pi )\bigl) .$
By the definitions, we have $\varphi _i (\pi )=\varepsilon _i (\pi )+
\alpha _i ^{\vee}\bigl( \mathrm{wt}(\pi )\bigl)$ for all $i\in I$.
Next, we define the Kashiwara operators.
We use the root operators $e_i ,\ i\in I^{re},$ and 
$f_i ,\ i\in I,$ on $\mathbb{P}$ as Kashiwara operators.
For $e_i ,\ i\in I^{im}$, we use the ``cutoff" of the
root operators $e_i ,\ i\in I^{im}$, on $\mathbb{P}\ \bigl( \supset \mathbb{B}(\lambda )\bigl)$,
i.e., if $e_i \pi \notin \mathbb{B}(\lambda )$, then we set $e_i \pi :=\mathbf{0}$ in 
$\mathbb{B}(\lambda )$ even if $e_i \pi \neq \mathbf{0}$ in $\mathbb{P}$.
Thus, $\mathbb{B}(\lambda )$ is endowed with a crystal structure. 
In particular, $\mathbb{B}(\lambda )$ is a \textit{normal\ crystal}
(for details, see [\textbf{JL}, \S 3.1.3]). 
Note that the crystal structure on $\mathbb{B}(\lambda )$
defined above agrees with that in [\textbf{JL}]. Thus, we know the following from [\textbf{JL}].
\begin{thm}$\mathrm{([\mathbf{JL},\ Proposition\ 6.3.5]).}$\ 
Let $\pi \in \mathbb{B}(\lambda )$. Then, 
$e_i \pi =\mathbf{0}$ in $\mathbb{B}(\lambda )$ for all $i\in I$
if and only if $\pi =\pi _{\lambda } :=(\lambda ;\ \! 0,1)$.
Moreover, we have $\mathbb{B}(\lambda )
=\mathcal{F}\pi _{\lambda } \setminus \{ \mathbf{0} \}$,
where $\mathcal{F}$ is the monoid generated by the Kashiwara operators $f_i $ for $i\in I$.
\label{3.3.6}
\end{thm}
For $\pi _1 ,\pi _2 ,\ldots ,\pi _n \in \mathbb{P}$,
define the \textit{concatenation of paths} 
$\pi =\pi _1 \otimes \pi _2 \otimes \cdots \otimes \pi _n \in \mathbb{P}$ 
by the formula: $\pi (t) := \sum _{l=1}^{m-1}\pi _l (1) + \pi _m (nt-m+1)$
for $\frac{m-1}{n} \le t \le \frac{m}{n}$ and $1\le m \le n.$
The following proposition follows immediately from the tensor 
product rule for crystals.
\begin{prop}$\mathrm{(cf.\ [\mathbf{Li1\textbf{-}2}]).}$\ 
Let $\lambda _1 ,\ldots ,\lambda _n \in P^+$. We set 
\begin{center}
$\mathbb{B}(\lambda _1) \otimes \cdots \otimes 
\mathbb{B}(\lambda _n) :=\{ \pi _1 \otimes 
\cdots \otimes \pi _n \in \mathbb{P}\ |\ 
\pi _m \in \mathbb{B}(\lambda _m )\ (m=1,\ldots ,n) \}$,
\end{center}
and define a crystal structure on 
$\mathbb{B}(\lambda _1) \otimes \cdots \otimes 
\mathbb{B}(\lambda _n) \ \bigl( \subset \mathbb{P} \bigl)$
in the same way as $\mathbb{B}(\lambda )$ above 
by using the cutoff of the root operators $e_i ,\ i\in I^{im}$. 
Then, the resulting crystal is isomorphic to the 
tensor product of the crystals $\mathbb{B}(\lambda _1),\ldots ,\mathbb{B}(\lambda _n)$.
\label{3.3.7}
\end{prop}
Let us show one important property of the operators $e_i,\ i\in I^{im}$, for later use. 
For this purpose, we need the following lemma, which is proved by the
same argument as for [\textbf{JL}, Lemma 5.3.5].
\begin{lem}
Let $\mu ,\nu \in \mathcal{W}\lambda $ be such that $\mu \ge \nu $, and $\alpha _i \in \varPi _{im}$.
If $\alpha _i ^{\vee}(\mu )=\alpha _i ^{\vee}(\nu )$, then $r_i ^{-1}\mu \ge r_i ^{-1}\nu $. 
Moreover, if there exists an $a$-chain for $(\mu ,\nu )$ for some rational number $a\in (0,1]$, 
then there exists an $a$-chain for $\bigl( r_i ^{-1}\mu ,r_i ^{-1}\nu \bigl)$.
\label{3.3.8}
\end{lem}
Note that the $a$-chain above for $(r_i ^{-1}\mu ,r_i ^{-1}\nu )$ is not
necessarily of shape $\lambda $.
\begin{cor}\label{3.3.9}
Fix $i\in I^{im}$, and let $\pi =Ff_i F' \pi _{\lambda }\in \mathbb{B}(\lambda )$,
where $F,F' \in \mathcal{F}$, and $F=f_{i_l} \cdots f_{i_2} f_{i_1}$ with
$f_{i_s} \neq f_i$ for $s=1,2,\ldots ,l$. Then, the following are equivalent:
\begin{itemize}
 \item[(1)]
$e_i \pi \neq \mathbf{0} \ \mathrm{in}\ \mathbb{B}(\lambda)$;
 \item[(2)]
$e_+ ^i (\pi )=f_- ^i (F' \pi _{\lambda })$;
 \item[(3)]
for each $s,\ 1\le s \le l$, either 
$\alpha _{i_s} ^{\vee}(\alpha _i)=0$, or
$$\bigl[0,f_- ^i (F'\pi _{\lambda }) \bigl)
\ \cap \ \bigl(f_+ ^{i_s}(f_{i_{s-1}}\cdots  f_{i_1} f_i F'\pi _{\lambda }), 
f_- ^{i_s}(f_{i_{s-1}} \cdots f_{i_1} f_i F' \pi _{\lambda }) \bigl) \ =\ \emptyset .$$
\end{itemize}
\end{cor}
$\mathbf{Proof.}$\ 
From Propositions \ref{3.3.2}, \ref{3.3.4}, and the fact 
that imaginary simple roots are anti-dominant integral weights (see [\textbf{JL}, Lemma 2.1.11]), 
it follows that $(2)$ and $(3)$ are equivalent and that $(1)$ implies $(2)$.
Hence it suffices to show that $(2)$ and $(3)$ imply $(1)$.

From $(2)$ and $(3)$, we have $e_i \pi \in \mathbb{P}$.
Now, we set $f_- ^i :=f_- ^i (F' \pi _{\lambda })$, 
$\pi =(\lambda _1 ,\ldots ,\lambda _q ,\ldots ,\lambda _k ;$
$a_0 ,\ldots ,f_- ^i =a_q ,\ldots ,a_k )$, and
$e_i \pi =(r_i ^{-1}\lambda _1 ,\ldots ,r_i ^{-1}\lambda _q , 
\lambda _{q+1} ,\ldots ,\lambda _k; \ \! a_0 ,\ldots ,f_- ^i =a_q ,\ldots ,a_k )$. 
Here, by (3), we have 
$\alpha _i ^{\vee}(\lambda _1) =\cdots =\alpha _i ^{\vee}(\lambda _q)$.
Therefore, by Lemma \ref{3.3.8}, there exists an $a_u$-chain
for $\bigl( r_i ^{-1} \lambda _u ,r_i ^{-1} \lambda _{u+1}\bigl)$
for each $u=1,2,\ldots ,q-1$. Hence it remains to show that there exists an 
$a_q$-chain for $(r_i ^{-1}\lambda _q ,\lambda _{q+1})$.
If we show this, it follows that $e_i \pi \in \mathbb{B}(\lambda )$.
Since $\mathrm{depth}_i ^0 (\lambda _q -\lambda _{q+1} ) >0$, 
there must be a positive root of the form $w\alpha _i $ for 
$w\in \mathcal{W}_{re}$ in the $a_q$-chain 
for $(\lambda _q ,\lambda _{q+1})$.
Also, we have $a_q \alpha _i ^{\vee}(\lambda _q )=1-a_{ii}$ by (3).
Consequently, by Lemma \ref{3.3.5}, there exists an $a_q$-chain for 
$(r_i ^{-1}\lambda _q ,\lambda _{q+1})$. \qed

\section{Embedding of path models and a crystal isomorphism theorem}
In this section, we introduce our main tool, which we call an
\textit{embedding of path models}. It provides an embedding of Joseph--Lamprou's path model 
for a generalized Kac--Moody algebra into Littelmann's path model 
for a certain Kac--Moody algebra constructed from the Borcherds--Cartan datum 
of the given generalized Kac--Moody algebra. By using this embedding, 
we give a new proof of the isomorphism theorem for path crystals, obtained in [\textbf{JL}, \S 7.4].
\subsection{Embedding of path models}
Let $\bigl( A=( a_{ij}) _{i,j\in I}, \varPi =\{ \alpha _i \} _{i\in I}, 
\varPi ^{\vee}=\{ \alpha _i ^{\vee} \} _{i\in I}, P,P^{\vee} \bigl)$ 
be a Borcherds--Cartan datum. From this, we construct a new \textit{Cartan datum}.
Recall that $\Tilde{I}=\{ (i,1) \} _{i\in I^{re}}\sqcup \{ (i,m) \} _{i\in I^{im}, m\in \mathbb{Z}_{\ge 1}}$, 
and define a Cartan matrix $\widetilde{A}:=(\Tilde{a}_{(i,m), (j,n)})_{(i,m), (j,n) \in \Tilde{I}}$ by
\begin{center}
$\begin{cases}
\Tilde{a}_{(i,m), (i,m)} :=2 
\ \mathrm{for}\ (i,m) \in \Tilde{I}, \\
\Tilde{a}_{(i,m), (j,n)} :=a_{ij} \ 
\mathrm{for}\ (i,m),(j,n)\in \Tilde{I},\ 
\mathrm{with}\ (i,m) \neq (j,n).
\end{cases}$
\end{center}
Let us denote by $\bigl( \widetilde{A}, \widetilde{\varPi }:=\{ \Tilde{\alpha} _{(i,m)} \} _{(i,m)\in \Tilde{I}} ,
\widetilde{\varPi }^{\vee}:=\{ \Tilde{\alpha} _{(i,m)} ^{\vee} \} _{(i,m)\in \Tilde{I}},
\widetilde{P},\widetilde{P}^{\vee} \bigl)$ the Cartan datum associated to $\widetilde{A}$,
where $\widetilde{\varPi }, \widetilde{\varPi }^{\vee}$ are the sets of simple roots and simple coroots, respectively, 
$\widetilde{P}$ is a weight lattice, and $\widetilde{P}^{\vee}$ is a coweight lattice.
Let $\Tilde{\mathfrak{g}}$ be the associated Kac--Moody algebra, 
$\Tilde{\mathfrak{h}}:=\widetilde{P}^{\vee}\otimes _{\mathbb{Z}} \mathbb{C}$
the Cartan subalgebra, and $\widetilde{W}$ the Weyl group.
Note that the Coxeter group $\mathfrak{W}$ (in Definition \ref{2.2.3}) 
is not necessarily isomorphic to the Weyl group $\widetilde{W}$.
\begin{lem}
If $A$ is symmetrizable, then so is $\widetilde{A}$.
\label{4.1.1}
\end{lem}
$\mathbf{Proof.}$\ 
Take a diagonal matrix $D=\mathrm{diag}(d_i )_{i\in I}$ such that $DA$ is symmetric. 
It follows that $d_i \Tilde{a}_{(i,m), (j,n)}=d_i a_{ij}= d_j a_{ji}=d_j \Tilde{a}_{(j,n), (i,m)}
\ \mathrm{for}\ (i,m), (j,n)\in \Tilde{I}$, with $(i,m) \neq (j,n)$.
Hence, if we set $\Tilde{d}_{(i,m)}:=d_i$ and 
$\widetilde{D}:=\mathrm{diag} \bigl( \Tilde{d}_{(i,m)} \bigl) _{(i,m)\in \Tilde{I}}$,
then $\widetilde{D}\widetilde{A}$ is symmetric. Thus, $\widetilde{A}$ is symmetrizable. 

\qed

\vspace{0.1in}
For each $\mu \in \mathcal{W}P^+ =\mathcal{W}_{re}P^+$, we take (and fix) an element
$\Tilde{\mu }\in \widetilde{P}$ such that
\begin{center}
$\Tilde{\alpha} _{(i,m)}^{\vee}(\Tilde{\mu })
=\alpha _i ^{\vee}(\mu )\ \mathrm{for\ all}\ (i,m) \in \Tilde{I}$.
\end{center}
Let $\widetilde{\mathbb{P}}$ denote the set of all piecewise-linear continuous 
maps $\Tilde{\pi }: [0,1]\longrightarrow (\Tilde{\mathfrak{h}}_{\mathbb{R}})^*$ 
such that $\Tilde{\pi }(0)=0$, and $\widetilde{\mathbb{B}}(\Tilde{\mu })$ 
the set of all (G)LS paths of shape $\Tilde{\mu }$
for $\Tilde{\mathfrak{g}}$. Let $\widetilde{\mathcal{F}}$
denote the monoid generated by the root operators 
$f_{(i,m)}, (i,m) \in \Tilde{I},$ and let
$\pi _{\lambda _1 ,\ldots ,\lambda _n }:=\pi _{\lambda _1}\otimes \cdots \otimes
\pi _{\lambda _n }\in \mathbb{P}\ \bigl( \mathrm{resp.,}\ 
\pi _{\Tilde{\lambda }_1 ,\ldots ,\Tilde{\lambda }_n }:=\pi _{\Tilde{\lambda }_1} \otimes \cdots \otimes
\pi _{\Tilde{\lambda }_n} \in \widetilde{\mathbb{P}} \bigl)$ 
for $\lambda _1 ,\ldots ,\lambda _n \in \mathcal{W}P^+$.
Following the notation of \S 2.2, we write $F_{\mathbf{i}}=f_{i_k }\cdots f_{i_2} f_{i_1} \in \mathcal{F}$ 
for $\mathbf{i}=( i_k ,\ldots ,i_2 ,i_1 ) \in \mathcal{I}$, and 
$F_{\mathbf{(i,m)}}=f_{(i_k , m_k )} \cdots f_{(i_2 , m_2)} f_{(i_1, m_1)}
\in \widetilde{\mathcal{F}}$ for 
$\mathbf{(i, m)}=((i_s , m_s)) _{s=1} ^k \in \widetilde{\mathcal{I}}$.
\begin{prop} For $\lambda _1 ,\ldots ,\lambda _n \in \mathcal{W}P^+$,
the map $\widetilde{\ \ }$ given below is well-defined and injective:
\begin{center}
$\widetilde{\ \ }:\ 
\mathcal{F}\pi _{\lambda _1 ,\ldots ,\lambda _n }
\longrightarrow \widetilde{\mathcal{F}}
\pi _{\Tilde{\lambda }_1 ,\ldots ,\Tilde{\lambda }_n },\ 
\pi = F_{\mathbf{i}}\pi _{\lambda _1 ,\ldots ,\lambda _n }
\longmapsto \Tilde{\pi }:= F_{\mathbf{(i,m)}}
\pi _{\Tilde{\lambda }_1 ,\ldots ,\Tilde{\lambda }_n }$,
\end{center}
where the $\mathbf{m}$ for which $\mathbf{(i,m)} \in \widetilde{\mathcal{I}}_{\mathrm{ord}}$ 
is determined uniquely by $\mathbf{i}\in \mathcal{I}$.
\label{4.1.2}
\end{prop}
The map $\widetilde{\ \ }$ in the proposition above is not 
necessarily a morphism of crystals. However, it is a 
\textit{quasi-embedding} of crystals in the sense that the equality
$\widetilde{f_i \pi }=f_{(i,m)} \Tilde{\pi }$ holds for
$\pi \in \mathcal{F}\pi _{\lambda _1 ,\ldots ,\lambda _n }$
if $i\in I^{re}$ and $m=1$, or for 
$\pi = F_{\mathbf{i}}\pi _{\lambda _1 , \ldots ,\lambda _n} 
\in \mathcal{F}\pi _{\lambda _1 ,\ldots ,\lambda _n }$ 
if  $i\in I^{im}$ and $m-1$ is the number of appearances of $i$ in $\mathbf{i}$.
\begin{ex}
Let 
$A=\left(
\begin{array}{@{\,}c|c@{\,}}
 2 & -1 \\ \hline
-2 & -4 \\ 
\end{array}\right)$
be a Borcherds--Cartan matrix, where 
$I=\{ 1, \mathbf{2} \}$, $I^{re}=\{ 1\}$, and $I^{im}=\{ \mathbf{2} \}$. We set
$\Tilde{I}=\{ (1,1) \} \sqcup \{ (\mathbf{2}, 1), (\mathbf{2}, 2), (\mathbf{2}, 3), (\mathbf{2}, 4), \ldots \}$, and define
\begin{center}
$\widetilde{A}=\left(
\begin{array}{@{\,}c|ccccc@{\,}}
 2 & -1 & -1 & -1 & -1 & \cdots \\ \hline
-2 &  2 & -4 & -4 & -4 & \cdots \\ 
-2 & -4 &  2 & -4 & -4 & \cdots \\ 
-2 & -4 & -4 &  2 & -4 & \cdots \\
-2 & -4 & -4 & -4 &  2 & \cdots \\
\vdots & \vdots & \vdots & \vdots & \vdots & \ddots 
\end{array}\right).$
\end{center}
Let 
$\pi =f_{\mathbf{2}} f_{\mathbf{2}} 
f_1 f_{\mathbf{2}} f_{\mathbf{2}} f_{\mathbf{2}} 
f_1 f_1 f_{\mathbf{2}} f_{\mathbf{2}} \pi _{\lambda }\in 
\mathbb{B}(\lambda )$ 
be a GLS path of shape $\lambda $. Then the corresponding LS path $\Tilde{\pi }$ 
of shape $\Tilde{\lambda }$ is given as follows:
\begin{center}
$\Tilde{\pi }=f_{(\mathbf{2}, 7)} f_{(\mathbf{2}, 6)} f_{(1,1)} 
f_{(\mathbf{2}, 5)} f_{(\mathbf{2},4)} f_{(\mathbf{2}, 3)} 
f_{(1,1)} f_{(1,1)} f_{(\mathbf{2},2)} f_{(\mathbf{2},1)} 
\pi _{\Tilde{\lambda }}\in \widetilde{\mathbb{B}}(\Tilde{\lambda} ).$
\end{center}
\label{4.1.3}
\end{ex}
Proposition \ref{4.1.2} is established by combining the following lemmas.

In what follows, we denote by $[\mu ]$ the unique element in 
$\mathcal{W}_{re} \mu \cap P^+$ for $\mu \in \mathcal{W}P^+$.
\begin{rem}
Note that $\pi _{\lambda }\in \mathbb{B}\bigl( [\lambda ] \bigl)$ 
$\bigl( \mathit{resp.,}\ \pi _{\Tilde{\lambda }} 
\in \widetilde{\mathbb{B}}\bigl( [\Tilde{\lambda }]\bigl) \bigl)$
for all $\lambda \in \mathcal{W}P^+ $. 
It follows that for each $\pi \in \mathcal{F}\pi _{\lambda _1 ,\ldots ,\lambda _n }$
$\bigl( \mathit{resp.,}\ \Tilde{\pi }\in \widetilde{\mathcal{F}}
\pi _{\Tilde{\lambda }_1 ,\ldots ,\Tilde{\lambda }_n } \bigl)$, 
we have $\pi \in \mathbb{B}\bigl( [\lambda _1]\bigl) 
\otimes \cdots \otimes \mathbb{B}\bigl( [\lambda _n ] \bigl)$
$\bigl( \mathit{resp.,}\ \Tilde{\pi } \in \widetilde{\mathbb{B}}\bigl( [\Tilde{\lambda} _1]\bigl) 
\otimes \cdots \otimes \widetilde{\mathbb{B}}\bigl( [\Tilde{\lambda} _n ] \bigl) \bigl)$.
Therefore, each of the tensor factors of these elements satisfies 
the chain condition (see Definition \ref{3.2.2}).
\label{4.1.4}
\end{rem}
Let $\mathbf{i}=(i_k ,\ldots ,i_2 ,i_1)\in \mathcal{I}$ and 
$\mathbf{m}=(m_k ,\ldots , m_2 ,m_1)$.
For each $s=1,2,\ldots ,k$, we set $\mathbf{i}_{[s]} :=(i_s ,\ldots ,i_2 ,i_1)\in \mathcal{I}$
and $\mathbf{(i,m)} _{[s]} :=((i_s , m_s), \ldots , (i_2 , m_2) , (i_1, m_1)) 
\in \widetilde{\mathcal{I}}$. Let us take an element
$\pi =F_{\mathbf{i}} \pi _{\Lambda } :=
F_{\mathbf{i}} \pi _{\lambda _1 ,\ldots ,\lambda _n }$,
and set $\Tilde{\pi }=F_{\mathbf{(i, m)}}
\pi _{\Tilde{\Lambda }} :=F_{\mathbf{(i,m)}} 
\pi _{\Tilde{\lambda }_1 ,\ldots ,\Tilde{\lambda }_n }$,
where $\Lambda := (\lambda _1 ,\ldots ,\lambda _n)$, 
$\Tilde{\Lambda} := (\Tilde{\lambda} _1 ,\ldots , \Tilde{\lambda} _n)$, and 
$\mathbf{(i,m)} \in \widetilde{\mathcal{I}}_{\mathrm{ord}}$.
In the following, we write $H[\pi ; i](t)$ instead of $H^{\pi } _i (t)$.
\begin{lem}\label{4.1.5}
With the notation above, if $\pi \in \mathbb{P} \setminus \{ \mathbf{0}\}$,
then $\Tilde{\pi }\in \widetilde{\mathbb{P}} \setminus \{ \mathbf{0}\}$. 
Moreover, we have $H[F_{\mathbf{i}_{[s-1]}} \pi _{\Lambda }; i_s](t)
=H[F_{\mathbf{(i, m)} _{[s-1]}}\pi _{\Tilde{\Lambda }}; (i_s , m_s)](t)$
for all $s=1,2, \ldots ,k,$ and $t\in [0,1]$.
\end{lem}
$\mathbf{Proof.}$\ 
Note that for a path $\pi \in \mathbb{P} \setminus \{ \mathbf{0}\}$, 
it depends only on the values $H[\pi ; i](1)$ and $m_i ^{\pi }$ 
whether $f_i \pi =\mathbf{0}$ or not (see Remark \ref{3.2.3} (2)).
We proceed by induction on $k$, and show the two assertions of the lemma simultaneously.
If $k=0$, the assertion is clear. So, we assume that $k\ge 1$.
Then we have $F_{\mathbf{i}_{[k-1]}}\pi _{\Lambda }\neq \mathbf{0}$, and hence,
by our induction hypothesis, $F_{\mathbf{(i,m)}_{[k-1]}}
\pi _{\Tilde{\Lambda }}\neq \mathbf{0}$ and 
$H[F_{\mathbf{i}_{[s-1]}} \pi _{\Lambda }; i_s ](t)
=H[F_{\mathbf{(i,m)} _{[s-1]}} \pi _{\Tilde{\Lambda }}; (i_s ,m_s)](t)$ for all $s=1,2,\ldots ,k-1,\ t\in [0,1]$.
Therefore, there exist functions
$\psi _s :[0,1]\longrightarrow [0,1]$ for $s=1,2,\ldots ,k-1$ such that 
$\bigl( F_{\mathbf{i}_{[k-1]}}\pi _{\Lambda }\bigl) (t)=\pi _{\Lambda }(t) -\sum _{s=1}^{k-1}
\psi _s (t) \alpha _{i_s}$ and $\bigl( F_{\mathbf{(i,m)} _{[k-1]}} 
\pi _{\Tilde{\Lambda }}\bigl) (t)=\pi _{\Tilde{\Lambda }}(t)
-\sum _{s=1}^{k-1}\psi _s (t) \Tilde{\alpha} _{(i_s , m_s)}$ for all $t\in [0,1]$ (see Remark \ref{3.1.1}).
Also, by the definitions, we have $\alpha _{i_k }^{\vee}(\lambda _l )
=\Tilde{\alpha} _{(i_k , m_k)} ^{\vee}( \Tilde{\lambda }_l )$
for all $l=1,\ldots ,n,$ and $\alpha _{i_k }^{\vee}(\alpha _{i_s})
=\Tilde{\alpha} _{(i_k , m_k )}^{\vee} (\Tilde{\alpha}_{(i_s , m_s)} \bigl)$
for all $s=1,\ldots ,k-1$. This implies that 
$H[F_{\mathbf{i}_{[k-1]}} \pi _{\Lambda }; i_k](t)
=H[F_{\mathbf{(i, m)} _{[k-1]}} \pi _{\Tilde{\Lambda }}; (i_k , m_k )](t)$ 
for all $t\in [0,1]$. From this, we deduce that 
$\Tilde{\pi }=F_{\mathbf{(i, m)}_{[k]}} \pi _{\Tilde{\Lambda }} \neq \mathbf{0}$. \qed
\begin{cor}
With the notation above, $\pi \in \mathbb{P} \setminus \{ \mathbf{0}\}$
if and only if $\Tilde{\pi } \in \widetilde{\mathbb{P}} \setminus \{ \mathbf{0}\}$.
\label{4.1.6}
\end{cor}

For $\mathbf{i}=(i_k ,\ldots ,i_2 ,i_1 )\in \mathcal{I}$,
we set $\mathbf{i}(i):=(x_p,\ldots ,x_2 , x_1 )$ if 
$\{ x_1, x_2 , \ldots , x_p \} =\{ x \ |\ i_x =i,\ 1\le x \le k \}$ and 
$1\le x_1 <x_2 <\cdots <x_p \le k$.
\begin{lem}
Let $\mathbf{i}=(i_k ,\ldots ,i_2 ,i_1),\ 
\mathbf{j}=(j_k ,\ldots ,j_2 ,j_1)\in I^k$.
Suppose that $F_{\mathbf{i}}\pi _{\Lambda }=F_{\mathbf{j}}\pi _{\Lambda }$
in $\mathcal{F}\pi _{\Lambda }$. If 
$\mathbf{i}(i)=(x_p ,\ldots ,x_2 , x_1)$ and 
$\mathbf{j}(i)=(y_p ,\ldots , y_2 ,y_1 )$ for $i\in I^{im}$, then we have
$f_{\pm}^i \bigl( F_{\mathbf{i}_{[x_q -1]}} 
\pi _{\Lambda } \bigl) =f_{\pm}^i \bigl( F_{\mathbf{j}_{[y_q -1]}}
\pi _{\Lambda } \bigl)$ for all $1\le q\le p$.
\label{4.1.7}
\end{lem}
$\mathbf{Proof.}$\ 
We set $f_{\pm} ^{i_{x_q}} = 
f_{\pm}^i \bigl( F_{\mathbf{i}_{[x_q -1]}} 
\pi _{\Lambda } \bigl)$ and $f_{\pm} ^{j_{y_q}} =
f_{\pm}^i \bigl( F_{\mathbf{j}_{[y_q -1]}}
\pi _{\Lambda } \bigl)$ for $1\le q\le p$.
Note that $\alpha _i \in \varPi _{im}$ is anti-dominant integral weight (see [\textbf{JL}, Lemma 2.1.11]).
Therefore, we have $f_{\pm}^{i_{x_p}} \le \cdots \le f_{\pm}^{i_{x_2}}\le f_{\pm}^{i_{x_1}},$ and 
$f_{\pm}^{i_{y_p}} \le \cdots \le f_{\pm}^{i_{y_2}}\le f_{\pm}^{i_{y_1}}$.
Also, we have $\pi _1 \otimes \cdots \otimes \pi _n =\pi ' _1 \otimes \cdots \otimes \pi ' _n$
if and only if $\pi _1 =\pi ' _1 ,\ldots ,\pi _n =\pi ' _n $.
Thus it is enough to check when $n=1$, for which the claim follows from 
Proposition \ref{3.3.2}. \qed
\begin{lem}\label{4.1.8}
If we write $\pi _1 =F_{\mathbf{i}}\pi _{\Lambda }$ and $\pi _2 =F_{\mathbf{j}}\pi _{\Lambda }$,
then $\pi _1 =\pi _2$ if and only if $\Tilde{\pi }_1 =\Tilde{\pi }_2$.
\end{lem}
$\mathbf{Proof.}$\ If $\pi _1 =\mathbf{0}$ or 
$\pi _2 =\mathbf{0}$, then the assertion follows from Corollary \ref{4.1.6}.
So, we may assume that $\pi _1 \neq \mathbf{0}$ and $\pi _2 \neq \mathbf{0}$.
Then, the same argument as in the proof of Lemma \ref{4.1.5} shows that
$\pi _1 (t)=\pi _{\Lambda }(t)-\sum _{s=1}^k \psi ^1 _s (t) \alpha _{i_s},$
$\pi _2 (t)=\pi _{\Lambda }(t)-\sum _{s=1}^k \psi ^2 _s (t) \alpha _{j_s},$
$\Tilde{\pi} _1 (t)=\pi _{\Tilde{\Lambda }}(t)-\sum _{s=1}^k \psi ^1 _s (t) \Tilde{\alpha} _{(i_s , m_s)},$ and 
$\Tilde{\pi} _2 (t)=\pi _{\Tilde{\Lambda }}(t)-\sum _{s=1}^k \psi ^2 _s (t) \Tilde{\alpha} _{(j_s , l_s)}$
for some functions $\psi _s ^{\epsilon }:[0,1]\longrightarrow [0,1] $, 
$s=1,\ldots ,k$ and $\epsilon =1,2$, where 
$\mathbf{i}=(i_k ,\ldots ,i_1) ,\ \mathbf{j}=(j_k , \ldots , j_1) \in \mathcal{I}$, and 
$\mathbf{(i, m)}=((i_s , m_s ))_{s=1} ^k, 
\mathbf{(j, l)}=((j_s , l_s ))_{s=1} ^k \in \widetilde{\mathcal{I}}_{\mathrm{ord}}$.
 
Suppose that $\pi _1 =\pi _2 $. Then we have
$$\sum _{1\le s\le k, i_s =i} \psi ^1 _s (t)=\sum _{1\le s\le k, j_s =i} \psi ^2 _s (t).$$
Therefore, for $i\in I^{re}$, we see that
$$\sum _{ \begin{subarray}\ \ \ \ 1\le s\le k \\
(i_s , m_s)=(i, 1)\end{subarray}}
\psi ^1 _s (t) =\sum _{\begin{subarray}\ \ 1\le s\le k \\
\ \ i_s =i \end{subarray}}
\psi ^1 _s (t) =\sum _{ \begin{subarray}\ \ 1\le s\le k \\
\ \ j_s =i \end{subarray}}
\psi ^2 _s (t) =\sum _{ \begin{subarray}\ \ \ \ 1\le s\le k \\
\ (j_s , l_s)=(i,1) \end{subarray}} \psi ^2 _s (t).$$
This implies that the coefficients of $\Tilde{\alpha} _{(i,1)}$ 
are equal for $\Tilde{\pi }_1 (t)$ and $\Tilde{\pi }_2 (t)$.

Next, let $i\in I^{im}$.
If we write $\mathbf{i}(i)=\bigl( x_p ,\ldots , x_1 \bigl)$
and $\mathbf{j}(i)=\bigl( y_p ,\ldots , y_1 \bigl)$,
then we have $\psi ^1 _{x_q} (t) =\psi ^2 _{y_q} (t)$
for each $1\le q\le p$ by Proposition \ref{3.3.2} and Lemma \ref{4.1.7}. 
Hence for each $1\le q \le p$, the coefficients of $\Tilde{\alpha } _{(i,q)}$ are equal
for $\Tilde{\pi }_1 (t)$ and $\Tilde{\pi }_2 (t)$.
This shows the equality $\Tilde{\pi }_1 =\Tilde{\pi }_2$.

Conversely, suppose that $\Tilde{\pi }_1 =\Tilde{\pi }_2$. 
We can show the equality $\pi _1 =\pi _2$
in a way similar to the above. This proves the lemma. \qed

\subsection{An isomorphism theorem}
Let $\lambda _1 ,\ldots ,\lambda _n \in \mathcal{W}P^+$.
If $\pi _{\Lambda }=\pi _{\lambda _1 ,\ldots ,\lambda _n }$
is a \textit{dominant path}, i.e., $\pi _{\Lambda }(t) \in \sum _{\eta \in P^+}\mathbb{R}_{\ge 0}\eta $
for all $t\in [0,1]$, then $\pi _{\Tilde{\Lambda }}=\pi _{\Tilde{\lambda }_1 ,\ldots ,\Tilde{\lambda }_n }$
is also a dominant path. In particular, it follows that $\lambda =\lambda _1 +\cdots +\lambda _n 
\ (\mathrm{resp.,}\ \Tilde{\lambda }= \Tilde{\lambda }_1 +\cdots +\Tilde{\lambda }_n)$ 
is a dominant integral weight.
\begin{thm}$\mathrm{(Isomorphism\ Theorem).}$\ 
With the notation above, we have the following isomorphism of crystals:
\begin{center} 
$\mathbb{B}(\lambda )
\xrightarrow{\ \cong \ } \mathcal{F}\pi _{\Lambda },\ 
F\pi _{\lambda }\longmapsto F\pi _{\Lambda }\ (F\in \mathcal{F})$.
\end{center}
\label{4.2.1}
\end{thm}
$\mathbf{Proof.}$\ 
Recall that $\mathbb{B}(\lambda )=\mathcal{F}\pi _{\lambda }$
and $\widetilde{\mathbb{B}}(\Tilde{\lambda })=\widetilde{\mathcal{F}}\pi _{\Tilde{\lambda }}$.
Hence, by Proposition \ref{4.1.2}, we have an embedding
$$\mathbb{B}(\lambda ) \hookrightarrow \widetilde{\mathbb{B}}(\Tilde{\lambda }),\ 
F_{\mathbf{i}}\pi _{\lambda }\mapsto F_{\mathbf{(i,m)}}\pi _{\Tilde{\lambda }}\ 
\bigl( \mathbf{(i,m)} \in \widetilde{\mathcal{I}}_{\mathrm{ord}} \bigl) .$$
Also, again by Proposition \ref{4.1.2}, we have an embedding
$$\mathcal{F}\pi _{\Lambda }\hookrightarrow \widetilde{\mathcal{F}}\pi _{\Tilde{\Lambda }},\ 
F_{\mathbf{i}}\pi _{\Lambda }\mapsto F_{\mathbf{(i,m)}}\pi _{\Tilde{\Lambda }}\ 
\bigl( \mathbf{(i,m)}\in \widetilde{\mathcal{I}}_{\mathrm{ord}} \bigl) .$$
Here, by the isomorphism theorem for path crystals for Kac--Moody algebras 
due to Littelmann ([\textbf{Li2},\ Theorem\ 7.1]), we have the isomorphism of crystals 
$$\widetilde{\mathbb{B}}(\Tilde{\lambda })\xrightarrow{\cong }\widetilde{\mathcal{F}}\pi _{\Tilde{\Lambda }},\ 
F \pi _{\Tilde{\lambda }}\mapsto F \pi _{\Tilde{\Lambda }}\ 
\bigl( F \in \widetilde{\mathcal{F}} \bigl) .$$
By composing these maps, we obtain the desired isomorphism. \qed

\section{A characterization of standard paths}
A path $\pi \in \mathbb{B}(\lambda _1 )\otimes 
\cdots \otimes \mathbb{B}(\lambda _n )$ is said to be \textit{standard}
if $\pi \in \mathcal{F}\pi _{\lambda _1 ,\ldots ,\lambda _n }$.
In this section, we give a necessary and sufficient
condition for a given path to be standard in terms of the monoid $\mathcal{W}$.
This result can be regarded as a generalization of [\textbf{Li3}, Theorem 10.1]
to the case of generalized Kac--Moody algebras. 

Let $\lambda _m \in P^+$, 
$\pi _m = (\lambda _{(m,1)}, \lambda _{(m,2)}, \ldots , \lambda _{(m, \kappa _m)}; \bm{a}_m ) 
\in \mathbb{B}(\lambda _m)$, 
$\bm{a}_m =(0<a_{(m,1)}< \cdots <a_{(m, \kappa _m)}=1)$, 
$\kappa _m \in \mathbb{Z}_{\ge 1}$ for $m=1,2,\ldots , n$, and set $\pi =\pi _1 \otimes \cdots \otimes \pi _n$.
Also, we define the set 
$\mathfrak{I}:=\{ (m,l) \mid l\in \{ 1,2,\ldots , \kappa _m \},\ m\in \{ 1,2,\ldots ,n \} \}
\subset  \mathbb{Z} \times \mathbb{Z}$
of indices of the weights appearing in the expressions for $\pi _m$, $1\le m \le n$, above, 
and equip the set $\mathfrak{I}$ with the lexicographic order, 
namely, we write $(m_1 , l_1) \ge (m_2 , l_2)$ if $m_1 > m_2$, or 
$m_1 = m_2 $ and $l_1 \ge l_2$. Throughout \S 5, we keep these notations.
\subsection{The connecting condition}
We introduce a certain condition for a concatenation of GLS paths, 
which we call the \textit{connecting\ condition}. Note that our condition below is similar to 
the one (the joining condition) introduced in [\textbf{JL}, \S 7.3.3], but is somewhat different. 
\begin{define}\label{5.1.1}
Let $\pi _m =(\lambda _{(m,1)}, \lambda _{(m,2)}, \ldots , \lambda _{(m, \kappa _m)}; \bm{a}_m) 
\in \mathbb{B}(\lambda _m )$, $m=1,2,\ldots ,n$, be as above. 
We say that the tuple $(\pi _1 ,\ldots ,\pi _n)$ satisfies the connecting condition if there exist
$w_{(m,l)}\in \mathcal{W}$ for $l=1,2,\ldots , \kappa _m $, $m=2,3,\ldots ,n-1 ,$ and $(m,l)=(n,1)$ such that
\begin{itemize}
\item[(1)] $w_{(m,l)} \lambda _m = \lambda _{(m,l)}$,  
\item[(2)] for each $1\le m \le n-1$,
$$\lambda _{(m, \kappa _m)}\ge w_{(m+1,1)}\lambda _m \ge w_{(m+1, 2)} \lambda _m \ge 
\cdots \ge w_{(m+1, \kappa _{m+1})} \lambda _m \ge w_{(m+2,1)} \lambda _m \ge \cdots$$
$$\cdots \ge w_{(n-1,1)} \lambda _m \ge w_{(n-1,2)} \lambda _m \ge 
\cdots \ge w_{(n-1, \kappa _{n-1})} \lambda _m \ge w_{(n,1)} \lambda _m \  
\text{in}\ \mathcal{W}\lambda _m,$$
\item[(3)] $w_{(m,1)} \lambda _m \in \mathcal{W}_{re}\lambda _m$ for each $1\le m \le n-1$,
\item[(4)] there exists a $1$-chain for 
$(\lambda _{(m, \kappa _m)}, w_{(m+1, 1)}\lambda _m)$ for each $1\le m \le n-1$.
\end{itemize}
\end{define}
Note that the elements $w_{(m+1,2)}\lambda _m , \ldots ,
w_{(n-1, \kappa _{n-1})} \lambda _m , w_{(n,1)}\lambda _m$ 
in the condition (2) above belong to
$\mathcal{W}_{re}\lambda _m$ for each $1\le m \le n-1$. We set
$$\mathcal{C}_{\lambda _1 ,\ldots ,\lambda _n }
:=\{ \pi _1 \otimes \cdots \otimes \pi _n \in \mathbb{B}(\lambda _1 )\otimes \cdots 
\otimes \mathbb{B}(\lambda _n ) \mid (\pi _1 ,\ldots ,\pi _n )\ \text{satisfies\ the\ connecting\ condition}\}.$$
The aim of this subsection is to prove the equality 
$\mathcal{C}_{\lambda _1 ,\ldots ,\lambda _n}
=\mathcal{F}\pi _{\lambda _1 ,\ldots ,\lambda _n }$.
For this purpose, we need to recall some fundamental properties 
of $a$-chains. For the proof of the following two lemmas, we refer the reader to [\textbf{JL, Li1-2}].
\begin{lem}$\mathrm{([\mathbf{JL,\ Li1\text{-}2}]).}$
Let $a\in (0,1]$ be a rational number, 
$\lambda \in P^+,\ i\in I^{re},$ and $\mu ,\nu \in \mathcal{W}\lambda $
with $\mu >\nu $. Then the following hold.
\begin{itemize}
\item[(1)]
If $r_i \mu <\mu$ and $r_i \nu \ge \nu$, then $\mu \ge r_i \nu \ (\ge \nu ).$
Moreover, if there exists an $a$-chain for $(\mu ,\nu )$, then 
there also exists one for $(\mu ,r_i \nu )$.
\item[(2)]
If $r_i \mu <\mu$ and $r_i \nu \ge \nu$, then $(\mu >)\ r_i \mu \ge \nu .$
Moreover, if there exists an $a$-chain for $(\mu ,\nu )$, then
there also exists one for $(r_i \mu ,\nu )$.
\item[(3)]
If $r_i \mu \le \mu$ and $r_i \nu > \nu$, then $(\mu \ge)\ r_i \mu \ge \nu.$
Moreover, if there exists an $a$-chain for $(\mu ,\nu )$, then
there also exists one for $(r_i \mu ,\nu )$.
\item[(4)]
If $r_i \mu \le \mu$ and $r_i \nu > \nu$, then $\mu \ge r_i \nu \ (>\nu ).$
Moreover, if there exists an $a$-chain for $(\mu ,\nu )$, then
there also exists one for $(\mu ,r_i \nu )$.
\item[(5)]
If $r_i \mu >\mu$ and $r_i \nu \ge \nu$, then $r_i \mu > r_i \nu \ (\ge \nu ).$
Moreover, if there exists an $a$-chain for $(\mu ,\nu )$, then
there also exists one for $(r_i \mu ,r_i \nu )$.
\item[(6)]
If $r_i \mu \le \mu$ and $r_i \nu < \nu$, then $(\mu \ge )\ r_i \mu > r_i \nu .$
Moreover, if there exists an $a$-chain for $(\mu ,\nu )$, then
there also exists one for $(r_i \mu ,r_i \nu )$.
\end{itemize}
\label{5.1.2}
\end{lem}
\begin{lem}$\mathrm{([\mathbf{JL},\ Lemma\ 7.3.1]).}$\ 
Let $\lambda ,\mu \in P^+$, $w\in \mathcal{W}$,
and $\beta \in \varDelta _{re}^+ \sqcup \varDelta _{im}$.
If $\beta ^{\vee}(w\lambda )>0$, then $\beta ^{\vee}(w\mu )\ge 0$.
\label{5.1.3}
\end{lem}
\begin{lem}
Let $\pi \in \mathcal{C}_{\lambda _1 ,\ldots ,\lambda _n }$ and $i\in I$.
If $f_i \pi \in \mathbb{P}$, then $f_i \pi \in \mathcal{C}_{\lambda _1 ,\ldots ,\lambda _n }$.
\label{5.1.4}
\end{lem}
$\mathbf{Proof.}$\ 
If $f_i \pi \in \mathbb{P}$, then $f_i \pi =
\pi _1 \otimes \cdots \otimes \pi _{m-1} \otimes 
(f_i \pi _m ) \otimes \pi _{m+1} \otimes \cdots 
\otimes \pi _n$ for some $m$, $1\le m\le n$, by Proposition \ref{3.3.7}.
\begin{flushleft}
\underline{Case $1$: $0<f_+ ^i (\pi _m ) <f_- ^i (\pi _m )<1.$}
\end{flushleft}
In this case, we necessarily have $i\in I^{re}$. Let us write $f_i \pi _m$ as
$(\lambda _{(m,1)}, \ldots , \lambda _{(m,p-1)}, r_i \lambda _{(m,p)}, \ldots ,$ 
$r_i \lambda _{(m,q)}, \lambda _{(m,q)}, \ldots , \lambda _{(m, \kappa _m)})$
with $1\le p \le q \le \kappa _m$. By conditions $(2)$ and $(3)$ in the connecting condition for $\pi $, 
we have $w_{(m,q)}\lambda _{m_1} \ge w_{(m,q+1)}\lambda _{m_1} \ge \cdots 
\ge w_{(m, \kappa _m)}\lambda _{m_1} \ge w_{(m+1,1)}\lambda _{m_1}$ in 
$\mathcal{W}_{re} \lambda _{m_1}$ for each $1\le m_1 <m$, 
and $\alpha _i ^{\vee}(w_{(m,q)}\lambda _m )=\alpha _i ^{\vee}(\lambda _{(m,q)})>0$ 
by the definition of root operators and the monotonicity of GLS paths (see Remark \ref{3.2.3} (2)).
Therefore, by Lemma \ref{5.1.3}, it follows that 
$\alpha _i ^{\vee}(w_{(m,q)}\lambda _{m_1})\ge 0$,
and hence $r_i w_{(m,q)}\lambda _{m_1} \ge w_{(m,q)}\lambda _{m_1}$
for each $1\le m_1 <m$. From this, we deduce that 
$r_i w_{(m,q)}\lambda _{m_1} \ge w_{(m,q)}\lambda _{m_1} \ge 
w_{(m,q+1)}\lambda _{m_1} \ge \cdots \ge w_{(m, \kappa _m)}\lambda _{m_1} 
\ge w_{(m+1,1)}\lambda _{m_1}$ in $\mathcal{W}_{re} \lambda _{m_1}$ for each $1\le m_1 <m$.
Hence, for $(m_1 , l) \ge (m,q)$ in $\mathfrak{I}$, we use $w_{(m_1, l)}$
given in the connecting condition for $\pi $.

Now, observe that by Remark \ref{3.2.3} (2), we have 
$\alpha _i ^{\vee}(w_{(m,l)}\lambda _m )
=\alpha _i ^{\vee}(\lambda _{(m,l)})>0$ for each $p\le l \le q$,
and hence $\alpha _i ^{\vee}(w_{(m,l)}\lambda _{m_1})\ge 0$
for $p\le l \le q$ and $1\le m_1 <m$ by Lemma \ref{5.1.3}.
Therefore, by Lemma \ref{5.1.2} (4) and (5), it follows that
$r_i w_{(m,p)}\lambda _{m_1} \ge r_i w_{(m,p+1)}\lambda _{m_1} \ge
\cdots \ge r_i w_{(m,q)}\lambda _{m_1}$ in $\mathcal{W}_{re}\lambda _{m_1}$
for each $m_1$, $1\le m_1 <m$.

We set 
$\mathfrak{J}_{(m, p)}:=\{ (m_1, l ) \in \mathfrak{I} \mid (m_1, l)<(m, p)\ \text{with}\ 
\alpha _i ^{\vee}(w_{(m_1, l)}\lambda _{m_2} )<0\ \text{for\ some}\ m_2,\ 1\le m_2 \le m_1 \}$. 
If $\mathfrak{J}_{(m, p)}$ is empty, then we have
$\alpha _i ^{\vee}(w_{(m_1, l)}\lambda _{m_2})\ge 0$
for all $(m_1, l)<(m, p)$ in $\mathfrak{I}$ and $1\le m_2 \le m_1$.
Therefore, by Lemma \ref{5.1.2} (4), (5) and (6), it follows that
$r_i \lambda _{(m_1, \kappa _{m_1})} \ge r_i w_{(m_1 +1,1)}\lambda _{m_1} \ge
\cdots \ge r_i w_{(m,p)}\lambda _{m_1}$ in $\mathcal{W}\lambda _{m_1}$ for each $1\le m_1 <m$.
In particular, there exists a $1$-chain for 
$(r_i \lambda _{(m_1, \kappa _{m_1})}, r_i w_{(m_1 +1,1)}\lambda _{m_1})$.
Also, in this case, the function $H_i ^{\pi }(t)$ is increasing in the interval 
$[a_{(1, \kappa_1 -1)}, f_+ ^i (\pi )]$.
However, the function $H_i ^{\pi }(t)$ attains its minimum at $f_+ ^i (\pi )$. 
From these, we deduce that it is constant in this interval. Therefore, we have
$\alpha _i ^{\vee}(w_{(m_1 , l)} \lambda _{m_1} )=0$ for all $(m_1, l)$, 
$(1, \kappa _1)\le (m_1, l)<(m, p)$ in $\mathfrak{I}$, and hence 
$r_i w_{(m_1, l)} \lambda _{m_1} =w_{(m_1, l)} \lambda _{m_1}$.
Thus, the connecting condition for $f_i \pi $ is satisfied by replacing 
$w_{(m_1, l)}$ with $r_i w_{(m_1, l)}$ for each $(1, \kappa _1)\le (m_1, l) \le (m, q)$.

Suppose that $\mathfrak{J}_{(m, p)}$ is nonempty, 
and let $(\Tilde{m}, \Tilde{l})$ be the maximum element in $\mathfrak{J}_{(m, p)}$.
Note that for $1\le m_1 <\Tilde{m}$, we have $r_i w_{(\Tilde{m}, \Tilde{l}+1)}\lambda _{m_1} \ge
\cdots \ge r_i w_{(m, p)}\lambda _{m_1}$ in $\mathcal{W}_{re}\lambda _{m_1}$, 
and $w_{(\Tilde{m}, \Tilde{l})}\lambda _{m_1} \ge r_i w_{(\Tilde{m}, \Tilde{l}+1)}\lambda _{m_1}$
by Lemma \ref{5.1.2} (1) and (4). Hence it follows that
$\lambda _{(m_1, \kappa _{m_1})} \ge w_{(m_1 +1,1)}\lambda _{m_1} 
\ge \cdots \ge w_{(\Tilde{m}, \Tilde{l})}\lambda _{m_1} \ge 
r_i w_{(\Tilde{m}, \Tilde{l}+1)}\lambda _{m_1} \ge \cdots \ge r_i w_{(m, p)}\lambda _{m_1}$
in $\mathcal{W}\lambda _{m_1}$.
Also, for $\Tilde{m} \le m_1 <m$, we have 
$r_i \lambda _{(m_1 , \kappa _{m_1})} \ge r_i w_{(m_1 +1,1)}\lambda _{m_1} \ge
\cdots \ge r_i w_{(m, p)} \lambda _{m_1}$ in $\mathcal{W}\lambda _{m_1}$.
In particular, if $(\Tilde{m}, \Tilde{l})=(m_1 , \kappa _{m_1})$, then there exists a $1$-chain for 
$(\lambda _{(m_1, \kappa _{m_1})}, r_i w_{(m_1 +1,1)}\lambda _{m_1})$ by Lemma \ref{5.1.2} (1) and (4).
Therefore, the same argument as above shows that 
$r_i w_{(m_1, l)}\lambda _{m_1}=w_{(m_1, l)}\lambda _{m_1}$
if  $(\Tilde{m}, \Tilde{l})<(m_1, l)<(m, p)$ in $\mathfrak{I}$.
Thus, if we replace each $w_{(m_1, l)}$ with $r_i w_{(m_1, l)}$ 
for $(\Tilde{m}, \Tilde{l}) <(m_1, l)<(m, q)$,
then the connecting condition for $f_i \pi $ is satisfied. 
\begin{flushleft}
\underline{Case $2$: $f_+ ^i (\pi _m )=0$.}
\end{flushleft}
If $i\in I^{re}$, then the assertion is shown in the same way as in Case 1.
So, we assume that $i\in I^{im}$. Note that $\alpha _i ^{\vee}(\lambda _{(m_1, \kappa _{m_1})})=0$ 
for $m_1<m$. Then, by condition $(2)$ in the connecting condition for $\pi $, we have
$0=\alpha _i ^{\vee}(\lambda _{(m_1, \kappa _{m_1})})
\ge \alpha _i ^{\vee}(w_{(m_1 +1,1)}\lambda _{m_1})\ge \cdots \ge 
\alpha _i ^{\vee}(w_{(n,1)}\lambda _{m_1})\ge 0$ for $m_1<m$, which implies that these are all equal to $0$,
and hence the elements $w_{(m_1+1,1)}\lambda _{m_1} , \ldots , w_{(n,1)}\lambda _{m_1}$
are fixed by $r_i$. In particular, $r_i w_{(m_1 +1,1)}\lambda _{m_1}=w_{(m_1+1,1)}\lambda _{m_1}$ 
belongs to $\mathcal{W}_{re}\lambda _{m_1}$. Thus, by replacing the (given) elements 
$w_{(m, l)}$ with $r_i w_{(m, l)}$ for each $1\le l \le q$, the connecting condition for $f_i \pi $ is satisfied.
\begin{flushleft}
\underline{Case $3$: $f_- ^i (\pi _m )=1$.}
\end{flushleft}
If we verify the existence of a $1$-chain for 
$(r_i \lambda _{(m, \kappa _m)}, w_{(m+1, 1)}\lambda _m)$,
then we can show that the connecting condition for $f_i \pi $ is satisfied 
in the same way as in Cases $1$ and $2$.
Since $f_- ^i (\pi _m )=1$, it follows that $\alpha _i ^{\vee}(\lambda _{(m, \kappa _m)})>0$, 
and that $r_i \lambda _{(m, \kappa _m)} \xleftarrow{\alpha _i } \lambda _{(m, \kappa _m)}$
is a $1$-chain. Therefore, we can combine a $1$-chain for
$(\lambda _{(m, \kappa _m)}, w_{(m+1,1)}\lambda _m)$ given in the 
connecting condition for $\pi $ with 
$r_i \lambda _{(m, \kappa _m)} \xleftarrow{\alpha _i }\lambda _{(m, \kappa _m)}$
to obtain a desired one.  \qed

\vspace{2mm}
The proof of the following lemma is similar to that of Lemma \ref{5.1.4}.
\begin{lem}\label{5.1.5}
Let $\pi \in \mathcal{C}_{\lambda _1 ,\ldots ,\lambda _n }$
and $i\in I^{re}$. If $e_i \pi \in \mathbb{P}$, then 
$e_i \pi \in \mathcal{C}_{\lambda _1 ,\ldots ,\lambda _n }$.
\end{lem}
\begin{lem}\label{5.1.6}
Let $\pi \in \mathcal{C}_{\lambda _1 ,\ldots ,\lambda _n }$
and $i\in I$. If $f_i \pi \ (\text{resp.,}\ e_i \pi ) \in 
\mathbb{B}(\lambda _1 )\otimes \cdots \otimes \mathbb{B}(\lambda _n )$,
then $f_i \pi \ (\text{resp.,}\ e_i \pi ) \in \mathcal{C}_{\lambda _1 ,\ldots ,\lambda _n }$.
\end{lem}
$\mathbf{Proof.}$\ 
Note that $\mathbb{B}(\lambda _1 )\otimes \cdots \otimes 
\mathbb{B}(\lambda _n )$ is closed in $\mathbb{P}$ 
by $e_i ,i\in I^{re}$, and $f_i ,i\in I$, i.e., if $e_i \pi \in \mathbb{P}$ for $i\in I^{re}$
(resp., $f_i \pi \in \mathbb{P}$ for $i\in I$), then $e_i \pi \ (\text{resp.,}\ f_i \pi ) \in 
\mathbb{B}(\lambda _1 )\otimes \cdots \otimes \mathbb{B}(\lambda _n )$ for 
$\pi \in \mathbb{B}(\lambda _1 )\otimes \cdots \otimes \mathbb{B}(\lambda _n )$.
From this, the assertion for $e_i ,i\in I^{re},$ and $f_i ,i\in I,$ follows.
Hence it suffices to show the assertion for $e_i ,\ i\in I^{im}$.

Let $i\in I^{im}$. By Propositions \ref{3.3.4} and \ref{3.3.7}, we can write
$e_i \pi =\pi _1 \otimes \cdots \otimes \pi _{m-1} \otimes 
(e_i \pi _m )\otimes \pi _{m+1} \otimes \cdots \otimes 
\pi _n \in \mathbb{B}(\lambda _1 )\otimes \cdots 
\otimes \mathbb{B}(\lambda _n )$,
where $e_i \pi _m =( r_i ^{-1} w_{(m,1)}\lambda _m , \ldots , 
r_i ^{-1} w_{(m,q )}\lambda _m , w_{(m, q+1)}\lambda _m ,$
$\ldots ,w_{(m, \kappa _m)}\lambda _m)$.
We will show that $e_i \pi $ satisfies the connecting condition.

Note that $\alpha _i ^{\vee}(\lambda _{(m_1, \kappa _{m_1})})=0$ for $1\le m_1 <m$. 
Then, by condition $(2)$ in the connecting condition for $\pi $, we have
$0=\alpha _i ^{\vee}(\lambda _{(m_1, \kappa _{m_1})})
\ge \alpha _i ^{\vee}(w_{(m_1 +1,1)}\lambda _{m_1})\ge
\cdots \ge \alpha _i ^{\vee}(w_{(n,1)}\lambda _{m_1})\ge 0$ for $1\le m_1 <m$, 
which implies that these are all equal to $0$, and hence 
$r_i ^{-1} w_{(m,l)} \lambda _{m_1}=w_{(m,l)}\lambda _{m_1}$ for all $1\le l\le q$.
In particular, we have $r_i ^{-1}w_{(m,1)}\lambda _{m-1}=w_{(m,1)}\lambda _{m-1} 
\in \mathcal{W}_{re}\lambda _{m-1}$ by condition $(3)$ in the connecting condition for $\pi $.

If $e_+ ^i (\pi _m )<1$, then the assertion is shown by replacing the (given) elements 
$w_{(m,l)}$, $l=1,2,\ldots ,q$, with the corresponding elements 
$r_i ^{-1}w_{(m,l)}$, $l=1,2,\ldots ,q$. Hence we assume that $e_+ ^i (\pi _m )=1$.
We must verify the existence of a $1$-chain for 
$(r_i ^{-1}\lambda _{(m, \kappa _m)}, w_{(m+1, 1)}\lambda _m)$.
By Proposition \ref{3.3.4}, we have 
$\alpha _i ^{\vee}(\lambda _{(m, \kappa _m)})=1-a_{ii}$.
Since $e_i \pi _m \in \mathbb{B}(\lambda _m)$,
it follows that $r_i ^{-1} \lambda _{(m, \kappa _m)} \in \lambda _m -Q^+$, and hence
$\mathrm{depth}_i ^{\lambda _m }(\lambda _{(m, \kappa _m)})>0$.
Also, by condition $(3)$ in the connecting condition, we have 
$w_{(m+1, 1)}\lambda _m \in \mathcal{W}_{re}\lambda _m 
\subset \lambda _m -\sum _{j\in I^{re}}\mathbb{Z}_{\ge 0}\alpha _j $.
Therefore, $\lambda _{(m, \kappa _m)} \neq w_{(m+1, 1)}\lambda _m$, and there exists a positive root 
of the form $w\alpha _i $ for $w\in \mathcal{W}_{re}$ in the $1$-chain for 
$(\lambda _{(m, \kappa _m)}, w_{(m+1,1)}\lambda _m)$, given by 
condition $(4)$ in the connecting condition for $\pi $. 
Consequently, by Lemma \ref{3.3.5}, there exists a $1$-chain for 
$(r_i ^{-1} \lambda _{(m, \kappa _m)}, w_{(m+1, 1)}\lambda _m)$.
Thus, in this case, if we replace the elements $w_{(m,l)}$, $l=1,2,\ldots , \kappa _m$,
which are given by the connecting condition for $\pi $, 
with the elements $r_i ^{-1}w_{(m,l)}$, $l=1,2,\ldots , \kappa _m$, 
then $e_i \pi $ satisfies connecting condition. \qed

\vspace{0.1in}
Let $\pi \in \mathbb{P}$ and $\lambda \in P$.
We say that $\pi $ is \textit{$\lambda $-dominant} if $\pi +\lambda $ is a dominant path,
where $(\pi +\lambda )(t):=\pi (t)+\lambda $ is a path in $\mathfrak{h}_{\mathbb{R}}^*$.
\begin{lem}\label{5.1.7}
Let $\pi =\pi _1 \otimes \cdots \otimes \pi _n 
\in \mathcal{C}_{\lambda _1 ,\ldots ,\lambda _n }$.
If $e_i \pi =\mathbf{0}$ in $\mathbb{B}(\lambda _1 )\otimes \cdots 
\otimes \mathbb{B}(\lambda _n )$ for all $i\in I$, 
then $\pi _m =\pi _{\lambda _m}$ for all $m=1,2,\ldots ,n$.
\end{lem}
$\mathbf{Proof.}$\ 
First we show that $\pi _1 =\pi _{\lambda _1}$.
Let $i\in I^{re}$. If $e_i \pi =\pi _1 \otimes e_i ( \pi _2 \otimes 
\cdots \otimes \pi _n )$, then $e_i \pi $ is not $\mathbf{0}$ by the 
normality of $\mathbb{B}(\lambda _1)\otimes \cdots \otimes \mathbb{B}(\lambda _n)$.
Therefore, we must have $e_i \pi =(e_i \pi _1 )\otimes \pi _2 \otimes \cdots \otimes \pi _n $.
In particular, it follows that
\begin{align}\label{eqA}
\varphi _i (\pi _1) \ge 
\varepsilon _i (\pi _2 \otimes \cdots \otimes 
\pi _n ). \tag{A}
\end{align}
This implies that $e_i \pi _1 =\mathbf{0}$ in $\mathbb{B}(\lambda _1)$.
Let $i\in I^{im}$. If $\alpha _i ^{\vee}\bigl( \mathrm{wt}(\pi _1)\bigl) \ge 
1-a_{ii}$, then $e_i \pi =(e_i \pi _1)\otimes \pi _2
\otimes \cdots \otimes \pi _n$ by the tensor product rule for crystals.
Therefore, we must have $e_i \pi _1 =\mathbf{0}$. 
If $\alpha _i ^{\vee}\bigl( \mathrm{wt}(\pi _1)\bigl) <1-a_{ii}$, then 
$e_i \pi _1 =\mathbf{0}$ by the definition of $e_i$.
Hence we have $e_i \pi _1 =\mathbf{0}$ in $\mathbb{B}(\lambda _1)$.
Consequently, we deduce from Theorem \ref{3.3.6} that $\pi _1 =\pi _{\lambda _1}$.
Moreover, by inequality (\ref{eqA}), it follows that
$H_j ^{\pi _2 \otimes \cdots \otimes \pi _n}(t)
+\alpha _j ^{\vee}(\lambda _1)\ge 
m_j ^{\pi _2 \otimes \cdots \otimes \pi _n}
+\alpha _j ^{\vee}(\lambda _1)
=-\varepsilon _j (\pi _2 \otimes \cdots 
\otimes \pi _n )+\varphi _j (\pi _{\lambda _1})\ge 0$
for all $j\in I^{re}$.
Thus, $\pi _2 \otimes \cdots \otimes \pi _n $ is $\lambda _1 $-dominant.

Now, we assume that $\pi =\pi _{\lambda _1}\otimes \cdots \otimes 
\pi _{\lambda _{m-1}}\otimes \pi _m \otimes \cdots \otimes \pi _n $,
and that $\pi _m \otimes \cdots \otimes \pi _n $
is $(\lambda _1 +\cdots +\lambda _{m-1})$-dominant.
We proceed by induction on $m$, and show that 
$\pi _m =\pi _{\lambda _m}$ and $\pi _{m+1}\otimes \cdots \otimes \pi _n $
is $(\lambda _1 +\cdots +\lambda _m )$-dominant.
We set $J:=\{ j\in I \mid \alpha _j ^{\vee}(\lambda _1 +\cdots +\lambda _{m-1})=0 \}$.
By condition $(2)$ in the connecting condition for $\pi $, 
we have $\lambda _{m_1} \ge w_{(m,1)}\lambda _{m_1}$ for each $1\le m_1 < m$.
Since $\lambda _{m_1}$ is the minimum element of $\mathcal{W}\lambda _{m_1}$
in the partial order defined in \S 3.2, we have $\lambda _{m_1}=w_{(m,1)}\lambda _{m_1}$ 
for each $1\le m_1 < m$. In particular, it follows that 
$w_{(m,1)}\in \mathrm{Stab}_{\mathcal{W}}(\lambda _1 +\cdots +\lambda _{m-1})
=\langle r_j \mid j\in J \rangle _{\mathrm{monoid}}$ by Lemma \ref{2.2.11}.
\begin{flushleft}
\underline{Case 1: $i\in J$.}
\end{flushleft} 
If $i\in I^{re}$, then we have
$H_i ^{\pi _m \otimes \cdots \otimes \pi _n }(t)
=H_i ^{\pi _m \otimes \cdots \otimes \pi _n }(t)
+\alpha _i ^{\vee}(\lambda _1 +\cdots +\lambda _{m-1})\ge 0$
for all $t\in [0,1]$ by the induction hypothesis.
In particular, we have 
$H_i ^{\pi _m}(t)= H_i ^{\pi _m \otimes \cdots \otimes \pi _n}
\bigl( \frac{t}{n-m+1}\bigl) \ge 0$
for all $t\in [0,1]$, and hence 
$e_i \pi _m =\mathbf{0}$ in $\mathbb{B}(\lambda _m )$.

If $i\in I^{im}$, then we have
$\mathbf{0}=e_i (\pi _{\lambda _1}\otimes \cdots \otimes 
\pi _{\lambda _{m-1}} \otimes \pi _m \otimes 
\cdots \otimes \pi _n )=\pi _{\lambda _1 }\otimes \cdots \otimes 
\pi _{\lambda _{m-1}}\otimes e_i (\pi _m \otimes \cdots \otimes \pi _n )$.
Therefore, $e_i (\pi _m \otimes \cdots \otimes \pi _n )
=\mathbf{0}$, and hence $e_i \pi _m =\mathbf{0}$ in $\mathbb{B}(\lambda _m )$.
\begin{flushleft}
\underline{Case 2: $i\notin J$.}
\end{flushleft}
In this case, $\lambda _{(m,1)}=w_{(m,1)}\lambda _m 
\in \lambda _m -\sum _{j\in J}\mathbb{Z}_{\ge 0}\alpha _j $.
Since $\lambda _{(m,1)}>\lambda _{(m,2)}> \cdots >\lambda _{(m, \kappa _m)}$
in $\mathcal{W}\lambda _m$ and $\lambda _{(m,1)} \prec \lambda _{(m,2)} 
\prec \cdots \prec \lambda _{(m, \kappa _m)}$, we have $\lambda _{(m,1)}, 
\lambda _{(m,2)}, \ldots , \lambda _{(m, \kappa _m)}\in \lambda _m -\sum _{j\in J}
\mathbb{Z}_{\ge 0}\alpha _j $.
This implies that $\mathrm{wt}(\pi _m )\in \lambda _m -\sum _{j\in J}
\mathbb{Z}_{\ge 0}\alpha _j$. So, we have $\mathrm{wt}(\pi _m )+\alpha _i
\notin \lambda _m -Q^+$ since $i\notin J$.
Therefore, it follows that $e_i \pi _m =\mathbf{0}$ in $\mathbb{B}(\lambda _m)$.

Consequently, we have $\pi _m = \pi _{\lambda _m }$,
and $e_j (\pi _{\lambda _1}\otimes \cdots \otimes 
\pi _{\lambda _m}\otimes \pi _{m+1}\otimes \cdots \otimes \pi _n )
=e_j (\pi _{\lambda _1} \otimes \cdots \otimes \pi _{\lambda _m })
\otimes (\pi _{m+1}\otimes \cdots \otimes \pi _n)$ for each $j\in I^{re}$
by the normality of $\mathbb{B}(\lambda _1)\otimes \cdots \otimes \mathbb{B}(\lambda _n)$.
In particular, $\varphi _j (\pi _{\lambda _1}\otimes \cdots \otimes \pi _{\lambda _m})
\ge \varepsilon _j (\pi _{m+1}\otimes \cdots \otimes \pi _n )$ for all $j\in I^{re}$.
From this, we deduce that $H_j ^{\pi _{m+1} \otimes \cdots \otimes \pi _n}(t)
+\alpha _j ^{\vee}(\lambda _1 +\cdots +\lambda _m)
\ge m_j ^{\pi _{m+1} \otimes \cdots \otimes \pi _n}
+\alpha _j ^{\vee}(\lambda _1 +\cdots +\lambda _m) 
=-\varepsilon _j (\pi _{m+1} \otimes \cdots \otimes \pi _n)
+\varphi _j (\pi _{\lambda _1}\otimes \cdots \otimes 
\pi _{\lambda _m}) \ge 0$ for all $t\in [0,1]$,
which means that $\pi _{m+1} \otimes \cdots \otimes \pi _n$
is $(\lambda _1 +\cdots +\lambda _m)$-dominant. \qed

\vspace{0.1in}
The following is a corollary of the proof of Lemma \ref{5.1.7} above.
\begin{cor}
Let $\lambda ,\mu \in P^+$, and 
$\pi _1 \otimes \pi _2 \in \mathbb{B}(\lambda )
\otimes \mathbb{B}(\mu )$. 
If $e_i (\pi _1 \otimes \pi _2 )=\mathbf{0}$
in $\mathbb{B}(\lambda )\otimes \mathbb{B}(\mu )$
for all $i\in I$, then
$\pi _1 =\pi _{\lambda }$,
and $\pi _2 $ is $\lambda $-dominant.
\label{5.1.8}
\end{cor}

By combining the lemmas above, we finally obtain the following.
\begin{prop}
The equality $\mathcal{C}_{\lambda _1 ,\ldots ,\lambda _n }
=\mathcal{F}\pi _{\lambda _1 ,\ldots ,\lambda _n}$ holds.
\label{5.1.9}
\end{prop}

\subsection{A characterization of standard paths in terms of $\mathcal{W}$}
In this subsection, we give a condition for standardness 
in terms of the Bruhat order on $\mathcal{W}$.
It can be regarded as a generalization of [\textbf{Li3}, Theorem 10.1]
to the case of generalized Kac--Moody algebras. 
Before doing this, we remark that the monoid $\mathcal{W}$ has the \textit{Lifting Property};
the proof is similar to that of [\textbf{Hu}, Proposition in \S 5.9].
\begin{lem}$\mathrm{(cf.\ [\mathbf{Hu},\ Proposition\ in\ \S 5.9]).}$\ 
Let $w' , w \in \mathcal{W},$ with $w' \le w,$ and let $i\in I^{re}$.
Then, either $r_i w' \le w$ or $r_i w' \le r_i w$ holds.
\label{5.2.1}
\end{lem}
\begin{thm}
Let $\pi _m =(\lambda _{(m,1)}, \lambda _{(m,2)}, \ldots , \lambda _{(m, \kappa _m)}; \bm{a}_m) 
\in \mathbb{B}(\lambda _m )$, $m=1,\ldots ,n$, and 
$\pi =\pi _1 \otimes \cdots \otimes \pi _n \in \mathbb{B}(\lambda _1)\otimes 
\cdots \otimes \mathbb{B}(\lambda _n)$ be as in \S 5.1. 
Then, $\pi $ is standard if and only if there exist 
$w_{(m,l)} \in \mathcal{W},\ (m,l) \in \mathfrak{I},$ such that
\begin{itemize}
\item[(1)] $w_{(m,l)} \lambda _m =\lambda _{(m,l)}$ for $(m,l) \in \mathfrak{I}$,
\item[(2)] $w_{(1,1)} \ge \cdots \ge w_{(1, \kappa _1)} \ge w_{(2,1)} \ge 
\cdots \cdots  \ge w_{(n-1, \kappa _{n-1})} \ge w_{(n,1)} \ge \cdots 
\ge w_{(n, \kappa _n)}$ in $\mathcal{W}$,
\item[(3)] $w_{(m+1, 1)}\lambda _m \in \mathcal{W}_{re}\lambda _m$ for $m=1,2,\ldots ,n-1$, 
\item[(4)] there exists a $1$-chain for $(\lambda _{(m, \kappa _m)}, w_{(m+1, 1)}\lambda _m)$ 
for $m=1,2,\ldots ,n-1$.
\label{5.2.2}
\end{itemize}
\end{thm}
We call the ordered collection 
$(w_{(1,1)}, w_{(1,2)}, \ldots , w_{(1, \kappa _1)}, w_{(2,1)}, w_{(2,2)}, 
\ldots , w_{(n-1, \kappa _{n-1})}, w_{(n,1)}, w_{(n,2)},$
$\ldots , w_{(n, \kappa _n)})$ in this theorem 
a \textit{defining\ chain} for $\pi $. Note that conditions (3) and (4) above are 
automatically satisfied in the case of ordinary Kac--Moody algebras.
\begin{flushleft}
\textit{\textbf{Proof\ of\ Theorem\ \ref{5.2.2}.}}
\end{flushleft}
Observe that the natural surjection $\mathcal{W}\rightarrow \mathcal{W}\lambda $
respects the Bruhat order on $\mathcal{W}$ and 
the partial order on $\mathcal{W}\lambda $ defined in \S 3.2.
Hence the condition above is clearly sufficient for standardness by Proposition \ref{5.1.9}.
So, we will prove that the condition above is also necessary for standardness.
Note that it suffices to show that the condition 
above is preserved by $f_i ,\ i\in I$, since the tuple $(1,\ldots ,1)$ 
is a defining chain for $\pi _{\lambda _1 ,\ldots ,\lambda _n}$.

We assume that
$f_i \pi =\pi _1 \otimes \cdots \otimes \pi _{m-1}
\otimes (f_i \pi _m )\otimes \pi _{m+1} \otimes \cdots \otimes \pi _n $, 
and denote by $\bm{\beta}_{(m_1, l)}$,  $1\le l \le \kappa _{m_1}$, $1\le m_1 \le n$, 
the ordered collection $(r_{\beta _1},\ldots ,r_{\beta _s})$
consisting of the reflections for the positive roots appearing in a fixed chain for 
$(w_{(m_1, l)} ,w_{(m_1, l+1)})$: 
$w_{(m_1, l)} \xleftarrow{\beta _1}\cdots \xleftarrow{\beta _s} w_{(m_1, l+1)}$ (in $\mathcal{W}$).
Also, we denote by $\bm{\gamma}_{m_1}$, $1 \le m_1 \le n$,  
a similar tuple $(r_{\gamma _1 },\ldots ,r_{\gamma _u})$
corresponding to a fixed $1$-chain for 
$(w_{(m_1, \kappa _{m_1})}\lambda _{m_1}, w_{(m_1 +1, 1)}\lambda _{m_1})$: 
$w_{(m_1, \kappa _{m_1})}\lambda _{m_1} 
\xleftarrow{\gamma _1}\cdots \xleftarrow{\gamma _u} w_{(m_1 +1,1)}\lambda _{m_1}$ 
(in $\mathcal{W}\lambda _{m_1}$). 
Note that $\bm{\beta}_{(m_1, l)}$, 
$1\le l \le \kappa _{m_1}$, $1\le m_1 \le n$, recover a defining chain for $\pi $.
\begin{flushleft}
\underline{Case 1: $i\in I^{im}.$} 
\end{flushleft}
Suppose that $f_i \pi _m = (r_i \lambda _{(m,1)}=r_i w_{(m,1)}\lambda _m, \ldots ,
r_i \lambda _{(m, q)}=r_i w_{(m, q)}\lambda _m, \lambda _{(m, q)}, \ldots , \lambda _{(m, \kappa _m)})$.
Then, we have $\alpha _i ^{\vee}(\lambda _{(m,1)}) = \cdots = \alpha _i ^{\vee}(\lambda _{(m,q)})>0$
by Proposition \ref{3.3.2}.
Since $\alpha _i ^{\vee}(w_{(m_1, l)} \lambda _{m_1})=0$
for all $(m_1, l)<(m, 1)$ in $\mathfrak{I}$,
$r_i $ commutes with the reflections in $\bm{\beta}_{(m_1, l)}$ for $(m_1, l)<(m, q)$.
Also, we have $w_{(m_1, \kappa _{m_1})} \lambda _{m_1} 
\ge w_{(m_1 +1, 1)}\lambda _{m_1}$ for each $1\le m_1 <m$. 
This implies that $0=\alpha _i ^{\vee}(w_{(m_1, \kappa _{m_1})}\lambda _{m_1})\ge 
\alpha _i ^{\vee}(w_{(m_1 +1, 1)}\lambda _{m_1})\ge 0$,
and hence $\alpha _i ^{\vee}(w_{(m_1, \kappa _{m_1})}\lambda _{m_1})
=\alpha _i ^{\vee}(w_{(m_1 +1, 1)}\lambda _{m_1})=0$ for each $1\le m_1 <m$.
Moreover, $r_i $ commutes with the reflections in $\bm{\gamma}_{m_1}$, 
and $r_i w_{(m_1 +1, 1)}\lambda _{m_1}=w_{(m_1 +1, 1)}\lambda _{m_1} 
\in \mathcal{W}_{re}\lambda _{m_1}$.
Therefore, there exists a $1$-chain for $(r_i w_{(m_1, \kappa _{m_1})}
\lambda _{m_1}, r_i w_{(m_1 +1, 1)}\lambda _{m_1})
=(w_{(m_1, \kappa _{m_1})}\lambda _{m_1}, w_{(m_1 +1, 1)}\lambda _{m_1})$
by condition $(4)$ above. Combining these, we obtain a desired 
defining chain for $f_i \pi $, which is constructed as follows:
$\bigl( \bm{\beta}_{(1,1)}, \bm{\beta}_{(1,2)}, \ldots , \bm{\beta}_{(m, q-1)}, 
\{ r_i \} \cup \bm{\beta}_{(m,q)}, \bm{\beta}_{(m, q+1)},$
$\ldots , \bm{\beta}_{(n, \kappa _n -1)}\bigl)$, 
where we set $\{ r_i \} \cup \bm{\beta}_{(m,q)}:=( r_i , r_{\beta _1 },\ldots , r_{\beta _s})$
for $\bm{\beta}_{(m,q)}=(r_{\beta _1 },\ldots ,r_{\beta _s})$.
\begin{flushleft}
\underline{Case 2: $i\in I^{re}$.}
\end{flushleft}
Suppose that $f_i \pi _m =(\lambda _{(m,1)}, \ldots , \lambda _{(m, p-1)},
r_i \lambda _{(m, p)}, \ldots , r_i \lambda _{(m, q)}, \lambda _{(m, q)}, \ldots , \lambda _{(m, \kappa _m)})$
for $1\le p \le q \le \kappa _m $. Here, by Remark \ref{3.2.3} (2), we have 
$\alpha _i ^{\vee}(\lambda _{(m, l)})>0$ for each $p\le l \le q$,
and hence $r_i w_{(m, l)} \xleftarrow{\alpha _i } w_{(m, l)}$,
$r_i w_{(m, l)} \ge r_i w_{(m, l+1)}$ for each $p\le l < q$ by Lemma \ref{5.2.1}.
Now, we set $\mathfrak{K}_{(m,p)}:=\{ (m_1, l) \in \mathfrak{I} \mid 
(m_1, l) < (m, p),\ w_{(m_1, l)} \xleftarrow{\alpha _i } r_i w_{(m_1, l)} \}$. 
Suppose first that $\mathfrak{K}_{(m,p)}$ is nonempty, and 
denote by $(\Tilde{m}, \Tilde{l})$ the maximum element in $\mathfrak{K}_{(m,p)}$. Then, 
\begin{align}\label{eqB}
\text{if}\ (\Tilde{m}, \Tilde{l})<(m_1, l)<(m, p)\ 
\mathrm{in}\ \mathfrak{I},\ \text{then}\ 
w_{(m_1 , l)} \xrightarrow{\alpha _i }r_i w_{(m_1 , l)}.
\tag{B}
\end{align}
Therefore, by Lemma \ref{5.2.1} again, we have
$r_i w_{(m_1, l)} \ge r_i w_{(m_1, l+1)}$ 
(and $r_i w_{(m_1, \kappa _{m_1})} \ge r_i w_{(m_1 +1, 1)}$) 
for all $(m_1, l)$, $(\Tilde{m}, \Tilde{l})<(m_1 , l)<(m, p)$ in $\mathfrak{I}$, 
and $w_{(m_1, l)} \xrightarrow{\alpha _i} r_i w_{(m_1, l)}$ 
since $\alpha _i ^{\vee}(w_{(m_1, l)} \lambda _{m_1})\ge 0$.
Hence it follows that the function $H_i ^{\pi }(t)$ is increasing in 
$[a_{(\Tilde{m}, \Tilde{l})}, f_+ ^i (\pi )=a_{(m, p-1)}]$.
However, the function $H_i ^{\pi }(t)$ attains its minimum at $t=f_+ ^i (\pi )$. 
Thus, it must be constant in this interval. Namely, we have
\begin{align}\label{eqC}
\alpha _i ^{\vee}(w_{(m_1, l)}\lambda _{m_1})=0\ 
\mathrm{for\ all}\ (m_1, l),\ (\Tilde{m}, \Tilde{l})<(m_1, l)<(m, p)\
\mathrm{in}\ \mathfrak{I}.
\tag{C}
\end{align}
This implies that $r_i w_{(m_1, l)} \lambda _{m_1} = \lambda _{(m_1, l)}$
for $(\Tilde{m}, \Tilde{l})<(m_1, l)<(m, p)$.
In addition, we have $w_{(\Tilde{m}, \Tilde{l})}> w_{(\Tilde{m}, \Tilde{l}+1)}$, 
$w_{(\Tilde{m}, \Tilde{l})} \xleftarrow{\alpha _i } r_i w_{(\Tilde{m}, \Tilde{l})}$ 
and $w_{(\Tilde{m}, \Tilde{l}+1)} \xrightarrow{\alpha _i }r_i w_{(\Tilde{m}, \Tilde{l}+1)}$
in $\mathcal{W}$. From Lemma \ref{5.2.1}, it follows that 
$r_i w_{(\Tilde{m}, \Tilde{l}+1)} \le w_{(\Tilde{m}, \Tilde{l})}$ or 
$r_i w_{(\Tilde{m}, \Tilde{l}+1)} \le r_i w_{(\Tilde{m}, \Tilde{l})}$,
and either cases, we have $r_i w_{(\Tilde{m}, \Tilde{l}+1)} \le w_{(\Tilde{m}, \Tilde{l})}$ 
since $r_i w_{(\Tilde{m}, \Tilde{l})} < w_{(\Tilde{m}, \Tilde{l})}$.

To prove that the sequence $(w_{(1,1)}, w_{(1,2)}, \ldots , w_{(\Tilde{m}, \Tilde{l})}, 
r_i w_{(\Tilde{m}, \Tilde{l} +1)}, \ldots , r_i w_{(m, q)}, w_{(m, q+1)}, \ldots ,$
$w_{(n, \kappa _n)})$ gives a defining chain for $f_i \pi $, 
we verify the existence of a $1$-chain.
By conditions $(3)$ and $(4)$ in Theorem \ref{5.2.2}, 
there exists a $1$-chain for 
$(w_{(m_1, \kappa _{m_1})} \lambda _{m_1} , w_{({m_1}+1, 1)} \lambda _{m_1})$,
and $w_{({m_1}+1,1)} \lambda _{m_1} \in \mathcal{W}_{re} \lambda _{m_1}$ 
for each $\Tilde{m} < m_1 < m$. 
Therefore, we have $r_i w_{(m_1 +1, 1)} \lambda _{m_1} \in \mathcal{W}_{re}\lambda _{m_1}$
and $w_{(m_1 +1, 1)} \xrightarrow{\alpha _i } r_i w_{(m_1 +1, 1)}$ by (\ref{eqB}) above.
This implies that $r_i w_{(m_1 +1, 1)} \lambda _{m_1} \ge w_{(m_1 +1, 1)}\lambda _{m_1}$
in $\mathcal{W}\lambda _{m_1}$.
Also, we have $r_i w_{(m_1, \kappa _{m_1})} \lambda _{m_1}
=w_{(m_1, \kappa _{m_1})} \lambda _{m_1}$ by (\ref{eqC}) above.
Now, by using Lemma \ref{5.1.2} (4), we can verify, by case-by-case calculations,
that there exists a $1$-chain for  
$(r_i w_{(m_1 , \kappa _{m_1})} \lambda _{m_1}=w_{(m_1, \kappa _{m_1})} \lambda _{m_1},$
$r_i w_{(m_1 +1, 1)} \lambda _{m_1})$.

Finally, if $\Tilde{l}=\kappa _{\Tilde{m}}$ for the $(\Tilde{m}, \Tilde{l})$, then 
$w_{(\Tilde{m}, \Tilde{l})} \lambda _{\Tilde{m}} 
\ge r_i w_{(\Tilde{m}, \Tilde{l})} \lambda _{\Tilde{m}}$
since $w_{(\Tilde{m}, \Tilde{l})} \xleftarrow{\alpha _i } r_i w_{(\Tilde{m}, \Tilde{l})}$, 
and $w_{(\Tilde{m}+1, 1)} \lambda _{\Tilde{m}} \le r_i w_{(\Tilde{m}+1, 1)} \lambda _{\Tilde{m}}
\in \mathcal{W}_{re}\lambda _{\Tilde{m}}$
since $w_{(\Tilde{m}+1, 1)} \xrightarrow{\alpha _i } r_i w_{(\Tilde{m}+1, 1)}$.
Now, by using Lemma \ref{5.1.2} (1) and (4), we can verify, by case-by-case calculations, 
that there exists a $1$-chain for 
$(w_{(\Tilde{m}, \Tilde{l})} \lambda _{\Tilde{m}}, r_i w_{(\Tilde{m}+1, 1)} \lambda _{\Tilde{m}})$. 

Note that if $\mathfrak{K}_{(m,p)}= \emptyset$, then (\ref{eqB}) and (\ref{eqC}) 
hold for all $(m_1, l) < (m,p)$ in $\mathfrak{I}$, and the assertion follows 
by the same argument as above. \qed

\section{Crystal isomorphism between $\mathbb{B}(\lambda )$ and $B(\lambda )$}
In this section, we assume that the Borcherds--Cartan datum is symmetrizable and even.
\subsection{Proof of the isomorphism $\mathbb{B}(\lambda ) \cong B(\lambda )$}
Let $B(\lambda )$ denote the crystal basis of the integrable 
highest weight $U_q (\mathfrak{g})$-module $V(\lambda )$
with highest weight $\lambda \in P^+$ (see [\textbf{JKK}]).
In [\textbf{JL}], Joseph--Lamprou defined the map
\begin{center}
$\psi _{\lambda ,\lambda +\mu }:\mathbb{B}(\lambda )\longrightarrow 
\mathbb{B}(\lambda ) \otimes \mathbb{B}(\mu ),\ 
\pi \longmapsto \pi \otimes \pi _{\mu }$,
\end{center}
for all $\lambda ,\mu \in P^+$.
We have $\pi \otimes \pi _{\mu }
\in \mathcal{F}(\pi _{\lambda }\otimes \pi _{\mu })$
by Theorem \ref{5.2.2}, and 
$\mathcal{F}(\pi _{\lambda }\otimes \pi _{\mu })
\cong \mathbb{B}(\lambda +\mu )$ by Theorem \ref{4.2.1}.
Therefore, the map $\psi _{\lambda ,\lambda +\mu }$
induces an embedding $\mathbb{B}(\lambda )
\hookrightarrow \mathbb{B}(\lambda + \mu )$,
which is a morphism of crystals, except for weights.
Thus, if we write $\lambda \ge \mu $ when $\lambda -\mu \in P^+$ for $\lambda ,\mu \in P^+$,
then $\bigl \{ \mathbb{B}(\lambda )\ (\lambda \in P^+);\ 
\psi _{\mu ,\nu}\ (\mu ,\nu \in P^+) \bigl \}$ forms an inductive system. 
It follows that the inductive limit $\mathbb{B}(\infty ):=\underrightarrow{\lim}\ \! \mathbb{B}(\lambda )$ 
is naturally endowed with a crystal structure. In fact, Joseph--Lamprou ([\textbf{JL},\ \S 8]) proved that 
$\mathbb{B}(\infty)$ is isomorphic to the crystal basis $B(\infty )$ 
of the negative part $U_q ^- (\mathfrak{g})$ of the quantized universal enveloping 
algebra $U_q (\mathfrak{g})$ (for the precise statement, see [\textbf{JKK},\ \textbf{JKKS}, \textbf{JL}]). 
By the definition of $\mathbb{B}(\infty)$ above, there exists a unique embedding 
$\Psi _{\lambda }:\mathbb{B}(\lambda )\hookrightarrow \mathbb{B}(\infty )\otimes 
\mathcal{T}_{\lambda }$ of crystals for each $\lambda \in P^+$. Also, we can verify that 
the map $\mathbb{B}(\lambda )\rightarrow \mathbb{B}(\infty )\otimes 
\mathcal{T}_{\lambda } \otimes C,\ \pi \mapsto \Psi _{\lambda }(\pi )\otimes c$,
induces an isomorphism $\mathbb{B}(\lambda )
\xrightarrow{\cong }\mathcal{F}(\pi _{\infty }\otimes t_{\lambda } 
\otimes c)$ of crystals, 
where $\pi _{\infty}$ is the highest weight element of $\mathbb{B}(\infty )$,
and $\mathcal{T}_{\lambda }:=\{ t_{\lambda } \} ,\ C:=\{ c \}$ are crystals
introduced in [\textbf{JKKS}]. 
Here we recall from [\textbf{JKKS},\ \S 5] that 
$B(\lambda )\cong \mathcal{F}(b_{\infty }\otimes t_{\lambda } 
\otimes c) \subset B(\infty) \otimes \mathcal{T}_{\lambda } \otimes C$ as crystals, 
where $b_{\infty} \in B(\infty)$ is the unique element of weight $0$ in $B(\infty )$.
Thus, we have proved the following theorem.
\begin{thm}\ 
For $\lambda \in P^+$, we have an isomorphism of crystals $\mathbb{B}(\lambda )\cong B(\lambda ).$
\label{6.1.1}
\end{thm}
As a consequence, we obtain the following embedding by Proposition \ref{4.1.2}:
\begin{align*}
B(\lambda ) \hookrightarrow \widetilde{B}(\Tilde{\lambda }),
\widetilde{F}_{\mathbf{i}} u_{\lambda } \mapsto 
\widetilde{F}_{\mathbf{(i,m)}} \Tilde{u}_{\Tilde{\lambda }}, 
\ \text{with}\ \mathbf{(i,m)} \in \widetilde{\mathcal{I}}_{\mathrm{ord}},
\end{align*}
where $\widetilde{B}(\Tilde{\lambda })$ denotes the crystal basis of 
the irreducible highest weight $U_q (\Tilde{\mathfrak{g}})$-module 
$\widetilde{V}(\Tilde{\lambda })$ of highest weight $\Tilde{\lambda }\in \widetilde{P}^+$,
$\widetilde{F}_{\mathbf{i}}:=\Tilde{f}_{i_k} \cdots \Tilde{f}_{i_2} \Tilde{f}_{i_1}$,
$\widetilde{F}_{\mathbf{(i,m)}}:=
\Tilde{f}_{(i_k , m_k)} \cdots \Tilde{f}_{(i_2 , m_2)} \Tilde{f}_{(i_1 , m_1)}$
are monomials of the Kashiwara operators, and $u_{\lambda } \in V(\lambda )$, 
$\Tilde{u}_{\Tilde{\lambda }} \in \widetilde{V}(\Tilde{\lambda })$ are the highest weight vectors. 
This embedding is not a morphism of crystals; however, 
it is a quasi-embedding of crystals in the sense of \S 4.1
since we have $B(\lambda) \cong \mathbb{B}(\lambda) = \mathcal{F}\pi _{\lambda}$ 
and $\widetilde{B}(\Tilde{\lambda}) \cong \widetilde{\mathbb{B}}(\Tilde{\lambda}) 
=\widetilde{\mathcal{F}} \pi _{\Tilde{\lambda}}$ for $\lambda \in P^+$ as crystals.

\section{Decomposition rules for crystals of GLS paths}
In this section, we assume that the Borcherds--Cartan datum is symmetrizable and even.
\subsection{Decomposition Rule for tensor products}
Let $\lambda , \mu \in P^+$, and let $\pi \in \mathbb{B}(\mu )$.
Then, $e_i (\pi _{\lambda } \otimes \pi )=\mathbf{0}$ in $\mathbb{B}(\lambda )
\otimes \mathbb{B}(\mu )$ for all $i\in I^{re}$ if and only if 
$\pi $ is $\lambda $-dominant. For $e_i $, $i\in I^{im}$, 
we have the following criterion for whether 
$e_i (\pi _{\lambda } \otimes \pi )=\mathbf{0}$ or not.
\begin{lem}\label{7.1.1}
Let $\lambda ,\mu \in P^+$, and let $\pi \in \mathbb{B}(\mu )$
be a $\lambda $-dominant path. Then, 
$e_i (\pi _{\lambda } \otimes \pi )=\mathbf{0}$
for all $i\in I^{im}$ if and only if 
$\Tilde{\pi }\in \widetilde{\mathbb{B}}(\Tilde{\mu })$
is $\Tilde{\lambda } $-dominant.
\end{lem}
$\mathbf{Proof.}$\ 
In the following, we write $m(\pi ;i)$ instead of $m_i ^{\pi }$. 
Let $\pi :=F_{\mathbf{i}}\pi _{\mu }$, 
$\Tilde{\pi }=F_{\mathbf{(i,m)}}\pi _{\Tilde{\mu }}$
for $\mathbf{i}=(i_k ,\ldots ,i_2 ,i_1 )\in \mathcal{I}$.
Since $\pi $ is $\lambda $-dominant, we have 
$H_{(i,1)} ^{\pi _{\Tilde{\lambda }} \otimes \Tilde{\pi }}(t)
=H_i ^{\pi _{\lambda } \otimes \pi }(t)\ge 0$ for all $i\in I^{re}$.
Also, if $(i,m)$ does not appear in $\mathbf{(i,m)}$,
then it is clear that $H_{(i,m)} ^{\pi _{\Tilde{\lambda }} \otimes \Tilde{\pi }}(t) \ge 0$.

Let us take $i \in I^{im}$ appearing in $\mathbf{i}$.
If $\mathbf{i}(i)=( x_p , \ldots , x_2 , x_1 )$, then it follows that
$$m \bigl( \Tilde{\pi } ; (i_{x_q }, m_{x_q}) \bigl) =
m \bigl( F_{\mathbf{(i, m)}} \pi _{\Tilde{\mu }}; (i_{x_q }, m_{x_q}) \bigl) \ge
m \bigl( F_{\mathbf{(i, m)} _{[x_q]} } \pi _{\Tilde{\mu }}; (i_{x_q }, m_{x_q}) \bigl) =-1$$
for each $q=1,2,\ldots ,p$, where 
$\mathbf{(i,m)}_{[s]} =((i_l , m_l)) _{l=1} ^s$ for 
$s=1,2,\ldots k$ as in \S 4.1. 

If $\alpha _i ^{\vee}(\lambda ) >0$, then 
$e_i (\pi _{\lambda }\otimes \pi )=(e_i \pi _{\lambda })\otimes \pi 
=\mathbf{0}\otimes \pi =\mathbf{0}$.
Also, since $\Tilde{\alpha} _{(i_{x_q} , m_{x_q} )} ^{\vee }(\Tilde{\lambda })
=\alpha _i ^{\vee}(\lambda ) \in \mathbb{Z}_{>0}$, we have 
$\Tilde{\alpha} _{(i_{x_q} , m_{x_q} )} ^{\vee }(\Tilde{\lambda })
+m\bigl( \Tilde{\pi }; (i_{x_q } , m_{x_q}) \bigl) \ge 0$
for each $q=1,2,\ldots ,p$.

If $\alpha _i ^{\vee}(\lambda )=0$, then $\Tilde{\alpha} _{(i_{x_q}, m_q )} ^{\vee }(\Tilde{\lambda })=0$
for each $q=1,2,\ldots ,p$, and $e_i (\pi _{\lambda } \otimes \pi )=\pi _{\lambda } \otimes (e_i \pi )$.
Hence we have $e_i (\pi _{\lambda } \otimes \pi )=\mathbf{0}$
if and only if $e_i \pi =\mathbf{0}$. Now, we set 
$f_{\pm } ^{i_s }:=f_{\pm } ^{i_s }(F_{\mathbf{i}_{[s-1]}} \pi _{\mu })$ and 
$f_{\pm } ^{(i_s , m_s)}:=f_{\pm } ^{(i_s , m_s)}
(F_{\mathbf{(i,m)}_{[s-1]}} \pi _{\Tilde{\mu }})$
for each $s=1,2,\ldots ,k$.
By Corollary \ref{3.3.9}, we have
$e_i \pi =\mathbf{0}$ if and only if
there exists $i_s $, with $x_p <s \le k$, such that
$\alpha _{i_s}^{\vee}(\alpha _{i_{x_p}}) \neq 0$ and 
$[0,f_- ^{i_{x_p}}) \cap  \bigl( f_+ ^{i_s} , f_- ^{i_s} \bigl) \neq \emptyset$.
Since $f_- ^{i_{x_p}}\le f_- ^{i_{x_{p-1}}} \le \cdots \le
f_- ^{i_{x_1}}$, this is equivalent to
$[0,f_- ^{i_{x_q}}) \cap \bigl( f_+ ^{i_s} , f_- ^{i_s} \bigl) 
\neq \emptyset$ for all $q=1,\ldots ,p$.
Furthermore, by using the embedding of Proposition \ref{4.1.2}, 
we infer that there exists $(i_s , m_s)$, with 
$x_p < s \le k$, such that 
$\Tilde{\alpha} _{(i_s , m_s)} ^{\vee}(\Tilde{\alpha} _{(i_{x_q}, m_{x_q})}) \neq 0$ and
$\bigl[ 0, f_- ^{(i_{x_q}, m_{x_q})} \bigl) \cap 
\bigl( f_+ ^{(i_s , m_s)} , f_- ^{(i_s , m_s)} \bigl) 
\neq \emptyset $ for all $q=1,\ldots ,p$.
Since the $(i_{x_q} , m_{x_q}), \ q=1,\ldots ,p ,$ are 
all distinct, these conditions are equivalent to 
$m\bigl( \Tilde{\pi }; (i_{x_q }, m_{x_q}) \bigl) >
m\bigl(F_{\mathbf{(i,m)}_{[x_q]}} \pi _{\Tilde{\mu }}; (i_{x_q }, m_{x_q}) \bigl)$
for all $q=1,2,\ldots ,p$. In this case, we have 
$m\bigl( \Tilde{\pi }; (i_{x_q }, m_{x_q}) \bigl) =0$ since 
$m\bigl( \Tilde{\pi }; (i_{x_q }, m_{x_q}) \bigl) \in \mathbb{Z}_{\le 0}$ and 
$m\bigl(F_{\mathbf{(i,m)}_{[x_q]}} \pi _{\Tilde{\mu }}; (i_{x_q }, m_{x_q}) \bigl) =-1$, 
$\Tilde{\alpha} _{(i_{x_q}, m_{x_q})} ^{\vee}(\Tilde{\lambda })=\alpha _i ^{\vee}(\lambda )=0$ 
for all $q=1,\ldots ,p,$ and hence $\Tilde{\alpha} _{(i_{x_q}, m_{x_q})} ^{\vee}(\Tilde{\lambda })+
m\bigl( \Tilde{\pi }; (i_{x_q }, m_{x_q}) \bigl) =0\ (\ge 0)$ for all $q=1,\ldots ,p$.
This proves the lemma. \qed

\vspace{2mm}
By Corollary \ref{5.1.8}, we can restate this lemma as follows:
\begin{prop}
Let $\pi _1 \in \mathbb{B}(\lambda )$ and 
$\pi _2 \in \mathbb{B}(\mu )$.
Then, $e_i (\pi _1 \otimes \pi _2 )=\mathbf{0}$
in $\mathbb{B}(\lambda )\otimes \mathbb{B}(\mu )$
for all $i\in I$ if and only if $\pi _1 =\pi _{\lambda }$,
and $\Tilde{\pi }_2 $ is $\Tilde{\lambda }$-dominant
in $\widetilde{\mathbb{B}}(\Tilde{\mu })$.
\label{7.1.2}
\end{prop}
From this proposition, by using Theorem \ref{6.1.1}, we obtain the 
following decomposition rule.
\begin{thm}
Let $\lambda ,\mu \in P^+$. Then, we have an isomorphism of crystals:
$$\mathbb{B}(\lambda )\otimes \mathbb{B}(\mu )\ \cong \ 
\bigsqcup _
{\begin{subarray}{c}
\pi \in \mathbb{B}(\mu ) \\ 
\Tilde{\pi }\ \!\!:\ \!\! \Tilde{\lambda }
\text{-}\mathrm{dominant}
\end{subarray}}
\mathbb{B}\bigl( \lambda +\pi (1) \bigl).$$
\label{7.1.3}
\end{thm}

\subsection{Branching Rule for restriction to Levi subalgebras}
Let $S \subset I$ be a subset. We set
$S^{re}:=S\cap I^{re}$ and $S^{im}:=S\cap I^{im}$.
Also, we set $\widetilde{S}:=\{ (i,1) \} _{i\in S^{re}}\sqcup 
\{ (i,m) \} _{i\in S^{im}, m\in \mathbb{Z}_{\ge 1}}$.
Let us denote by $\mathfrak{g}_S$ (resp., 
$\Tilde{\mathfrak{g}}_{\widetilde{S}}$)
the Levi subalgebra of ${\mathfrak{g}}$ 
corresponding to $S$ (resp., the Levi subalgebra of $\mathfrak{\Tilde{g}}$ 
corresponding to $\widetilde{S}$).
We say that a path $\pi $ is $\mathfrak{g}_S$-dominant 
(resp., $\mathfrak{\Tilde{g}}_{\widetilde{S}}$-dominant)
if $\pi $ is a dominant path for $\mathfrak{g}_S$
(resp., $\mathfrak{\Tilde{g}}_{\widetilde{S}}$).
The proof of the following lemma is similar to that
of Lemma \ref{7.1.1}.
\begin{lem}
Let $\pi \in \mathbb{B}(\lambda )$. Then, 
$e_i \pi =\mathbf{0}$ in $\mathbb{B}(\lambda )$
for all $i\in S$ if and only if $\Tilde{\pi } \in 
\widetilde{\mathbb{B}}(\Tilde{\lambda })$ is 
$\mathfrak{\Tilde{g}}_{\widetilde{S}}$-dominant.
\label{7.2.1}
\end{lem}

Let us denote by $\mathbb{B}_S (\lambda )$ the set of 
all GLS paths of shape $\lambda $ for $\mathfrak{g}_S$.
From Lemma \ref{7.2.1}, by using Theorem \ref{6.1.1}, we obtain the 
following branching rule.
\begin{thm}
Let $\lambda \in P^+$. Then, we have an isomorphism of $\mathfrak{g}_S$-crystals:
$$\mathbb{B}(\lambda )\ \cong \ \bigsqcup _
{\begin{subarray}{c}\pi \in \mathbb{B}(\lambda) \\ 
\Tilde{\pi }\ \!\! :\ \!\! \Tilde{\mathfrak{g}}_{\widetilde{S}}
\text{-}\mathrm{dominant}
\end{subarray}} \mathbb{B}_S \bigl( \pi (1) \bigl) .$$
\label{7.2.2}
\end{thm}

\section{Appendix}
In this appendix, we give the (postponed) proofs of results stated in \S 2.2.
\begin{flushleft}
\textit{\textbf{Proof\ of\ Lemma\ \ref{2.2.5}.}}
\end{flushleft}
It is clear that the map $\sigma $ is well-defined.
To prove that it is bijective, we will construct its inverse.
Now we define
$$\xi :\mathfrak{V}\longrightarrow \mathcal{W},\ 
S_{\mathbf{(i,m)}} \longmapsto 
R_{\mathbf{i}},$$
where $\mathbf{(i,m)} \in \widetilde{\mathcal{I}}_{\mathrm{ord}}$.
If $\xi $ is well-defined, then 
it clearly provides the inverse of $\sigma $.

To verify the well-definedness of $\xi $,
we define the free group $\widetilde{\mathfrak{W}}$
generated by the symbols $\Tilde{s}_{(i,m)},\ (i,m) \in \Tilde{I}$, and then a map 
$$\Tilde{\xi} :\widetilde{\mathfrak{W}}\longrightarrow \mathcal{W},\ 
\widetilde{S}_{\mathbf{(i,m)}} \longmapsto 
R_{\mathbf{i}},$$
where $\widetilde{S}_{\mathbf{(i,m)}} 
:=\Tilde{s}_{(i_k , m_k)} \cdots \Tilde{s}_{(i_2 , m_2)}
\Tilde{s}_{(i_1 , m_1)}\in \widetilde{\mathfrak{W}}$
if $\mathbf{i}=(i_k ,\ldots ,i_2 ,i_1 )$ and
$\mathbf{m}=(m_k ,\ldots ,m_2 ,m_1 )$.
Obviously, $\Tilde{\xi }$ is well-defined.
Moreover, if $\mathbf{(i,m)} \in \widetilde{\mathcal{I}}$ is such that
\begin{align}
m_x \neq m_y \ \mathrm{for}\ i_x =i_y \in I^{im},\ \mathrm{with}\ x\neq y,\ 1\le x,y \le k,
\label{eq1}
\end{align}
then $\Tilde{\xi }$ commutes with 
\textit{nil-moves} and \textit{braid-moves} for
$\widetilde{S}_{\mathbf{(i,m)}}$, that is, if 
$\widetilde{S}_{\mathbf{(i,m)}}=\Tilde{v}
\Tilde{s}_{(i,m)}\Tilde{s}_{(i,m)} \Tilde{u}$, then
$$\Tilde{\xi }( \Tilde{v}\Tilde{s}_{(i,m)}\Tilde{s}_{(i,m)}\Tilde{u})
=\Tilde{\xi }(\Tilde{v} \Tilde{u}),$$
and if 
$\widetilde{S}_{\mathbf{(i,m)}}
=\Tilde{v}\underbrace{\Tilde{s}_{(i,m)}
\Tilde{s}_{(j,n)}\Tilde{s}_{(i,m)} \cdots }
_{l\ \mathrm{letters}}\Tilde{u}$ with $\mathrm{x}_{(i,m), (j,n)}=l$, then
$$\Tilde{\xi }\bigl( 
\Tilde{v}\underbrace{\Tilde{s}_{(i,m)}
\Tilde{s}_{(j,n)} \Tilde{s}_{(i,m)} \cdots }
_{l\ \mathrm{letters}}\Tilde{u} \bigl)
=\Tilde{\xi }\bigl( \Tilde{v}\underbrace{\Tilde{s}_{(j,n)}
\Tilde{s}_{(i,m)}\Tilde{s}_{(j,n)} \cdots }
_{l\ \mathrm{letters}}\Tilde{u} \bigl).$$
Note that if $\mathbf{(i,m)} \in \widetilde{\mathcal{I}}$
satisfies condition $(\ref{eq1})$, then
the number of appearances of $\Tilde{s}_{(i,m)}$, with $i\in I^{im}$, 
in the word $\widetilde{S}_{\mathbf{(i,m)}}$ for $S_{\mathbf{(i,m)}}$ 
is constant for each $m\in \mathbb{Z}_{\ge 1}$
during the process of nil-moves and braid-moves.
 In particular, it is at most one.
Hence condition $(\ref{eq1})$ is preserved by nil-moves and braid-moves.

Now, let $S_{\mathbf{(i,m)}}=S_{\mathbf{(j,n)}}$ be 
two expressions of a given element in $\mathfrak{V}$, with 
$\mathbf{(i,m)}, \mathbf{(j,n)} \in \widetilde{\mathcal{I}}_{\mathrm{ord}}$.
If we take the words $\widetilde{S}_{\mathbf{(i,m)}}$, 
$\widetilde{S}_{\mathbf{(j,n)}} \in \widetilde{\mathfrak{W}}$
for $S_{\mathbf{(i,m)}}, S_{\mathbf{(j,n)}}$, respectively, then 
$\Tilde{\xi }( \widetilde{S} _{\mathbf{(i,m)}})= R_{\mathbf{i}}$, 
$\Tilde{\xi }( \widetilde{S} _{\mathbf{(j,n)}}) = R_{\mathbf{j}}$ by the definitions.
Also, by the \textit{Word Property} (see [\textbf{BB}, Theorem 3.3.1])
of the Coxeter group $\mathfrak{W}$, every given expression can be transformed 
into a reduced expression by repeated application of nil-moves and braid-moves.
Now, we lift these moves for $S_{\mathbf{(i,m)}}$
and $S_{\mathbf{(j,n)}}$ to $\widetilde{\mathfrak{W}}$ as follows:
$\widetilde{S}_{\mathbf{(i,m)}} =
\Tilde{v}_0 \rightarrow \Tilde{v}_1 \rightarrow 
\cdots \rightarrow \Tilde{v}_p ,\ 
\widetilde{S}_{\mathbf{(j,n)}} =
\Tilde{u}_0 \rightarrow \Tilde{u}_1 \rightarrow 
\cdots \rightarrow \Tilde{u}_q ,$
where $\Tilde{v}_p$ and $\Tilde{u}_q$
are the corresponding reduced words, respectively.
Note that $\mathbf{(i,m)}, \mathbf{(j,n)} \in \widetilde{\mathcal{I}}_{\mathrm{ord}}$ 
satisfy condition (1). Then, from the arguments above, we deduce that
$\Tilde{\xi }( \widetilde{S}_{\mathbf{(i,m)}})
=\Tilde{\xi }( \Tilde{v}_0 ) 
=\Tilde{\xi }( \Tilde{v}_1 ) =\cdots
=\Tilde{\xi }( \Tilde{v}_p )$ and 
$\Tilde{\xi }( \widetilde{S}_{\mathbf{(j,n)}})
=\Tilde{\xi }( \Tilde{u}_0 ) 
=\Tilde{\xi }( \Tilde{u}_1 ) =\cdots 
=\Tilde{\xi }( \Tilde{u}_q ) .$
Again, by the Word Property, $\Tilde{v}_p$ and
$\Tilde{u}_q$ are transformed into each other by 
repeated application of braid-moves only. Therefore, the same argument
as above yields $\Tilde{\xi }( \Tilde{v}_p )=\Tilde{\xi }( \Tilde{u}_q )$.
Thus, we have $\Tilde{\xi }( \widetilde{S}_{\mathbf{(i,m)}})
=\Tilde{\xi }( \widetilde{S}_{\mathbf{(j,n)}})$, and hence $R_{\mathbf{i}}=R_{\mathbf{j}}$. 
This proves the lemma. \qed

\vspace{2mm}
Corollary \ref{2.2.7} is clear by Lemma \ref{2.2.5}.
We give the proof of Corollary \ref{2.2.8}.
Let $\{ \gamma _{(i,m)} \} _{(i,m)\in \Tilde{I}}$ be the 
set of simple roots associated with the Coxeter group $\mathfrak{W}$.
\begin{flushleft}
\textit{\textbf{Proof\ of\ Corollary\ \ref{2.2.8}.}}
\end{flushleft}
Let $v= r_{i_k} \cdots r_{i_1}$ and
$v' =r_{j_l}\cdots r_{j_1}$.
The condition $v=r_{\beta }v' =
wr_i w^{-1} r_{j_l}\cdots r_{j_1}$
implies that $\# \{ 1\le y\le l \mid j_y =i\}=p-1$.
Therefore, we have 
$\sigma (v)=\sigma (r_{\beta } v')
=\sigma (w) s_{(i, x_p)} \sigma (w^{-1})
\sigma (v')=s_{\gamma } \sigma (v')$,
where we set $\gamma = \sigma (w)(\gamma _{(i, x_p)})$
and $s_{\gamma }$ denotes the reflection with respect to 
the root $\gamma$. Also, we have 
$\ell _{\mathfrak{W}} \bigl( s_{\gamma } \sigma (v) \bigl)
=\ell _{\mathfrak{W}} \bigl( \sigma (v') \bigl)
=\ell _{\mathcal{W}} (v')<\ell _{\mathcal{W}} (v)
=\ell _{\mathfrak{W}} \bigl( \sigma (v) \bigl)$, and hence
\begin{align}
s_{\gamma } \sigma (v)
=s_{(i_k , m_k )}\cdots 
\widehat{s_{(i_u , m_u)}} \cdots s_{(i_1 , m_1)}
\label{eq2}
\end{align}
for some $(i_u , m_u) \in \Tilde{I}$, $1\le u \le k$, by the Strong Exchange Property of $\mathfrak{W}$.
In particular, $s_{(i, x_p)}$ appears at most once on the right-hand side of (\ref{eq2})
since $m_x ,\ x\in \{ 1\le x \le k \mid i_x =i \}$ are all distinct. 
However, we have $s_{\gamma } \sigma (v) =\sigma (v') \in \mathfrak{V}$, 
and $s_{(i, x_p)}$ does not appear in the expression 
$\sigma (v')=s_{(j_l , n_l)} \cdots s_{(j_1 , n_1)}$.
Consequently, by the Word Property of $\mathfrak{W}$, $s_{(i, x_p)}$ 
does not appear on the right-hand side of (\ref{eq2}). 
From this, we deduce that $s_{(i_u , m_u)} = s_{(i, x_p)}$.
Thus, the expression of an element on the right-hand side of (\ref{eq2}) satisfies the condition for $\mathfrak{V}$. 
By applying $\sigma ^{-1}$ to both sides of (\ref{eq2}), we obtain
$v' =\sigma ^{-1} \bigl( s_{\gamma } \sigma (v) \bigl)
=r_{i_k }\cdots \widehat{r_{i_{x_p}}} \cdots r_{i_1}$. \qed
\begin{flushleft}
\textit{\textbf{Proof\ of\ Lemma\ \ref{2.2.11}.}}
\end{flushleft}
The inclusion $\supset $ is clear. Hence we show the 
opposite inclusion $\subset $.
Take $w\in \mathrm{Stab}_{\mathcal{W}}(\lambda )$
with dominant reduced expression
$w=w_k r_{i_k} \cdots w_1 r_{i_1} w_0$, where
$i_1 ,\ldots ,i_k \in I^{im}$, and 
$w_0 ,\ldots ,w_k \in \mathcal{W}_{re}$.
We proceed by induction on $k$.
If $k=0$, it is clear by a (well-known) property of ordinary Coxeter groups.
Hence we assume that $k\ge 1$. In this case, $\mu :=r_{i_k} w_{k-1} \cdots r_{i_1} w_0 \lambda $
is dominant, and $w_k \mu =\lambda \in P^+$, which implies that 
$\mu =\lambda $ and $w_k \in \langle r_i \mid \alpha _i ^{\vee}(\lambda )=0 \rangle$.
Thus, we may assume that $w_k =1$.
Note that $\alpha _{i_s}^{\vee}(w_{s-1} r_{i_{s-1}} \cdots w_1 r_{i_1} w_0 \lambda )\ge 0$ 
for all $s=1,2,\ldots ,k$.
Here, by the assumption that $w\lambda =\lambda $ and the linear independence
of the simple roots, all of these inequalities must be equalities. In particular, we have
$\alpha _{i_k}^{\vee} (w_{k-1} r_{i_{k-1}} \cdots w_1 r_{i_1}w_0 \lambda )=0$.
Consequently, we deduce that $\lambda =w\lambda 
=r_{i_k} w_{k-1} r_{i_{k-1}} \cdots w_1 r_{i_1 }w_0 \lambda 
=w_{k-1} r_{i_{k-1}} \cdots w_1 r_{i_1} w_0 \lambda$.
Now, by induction, the desired inclusion is shown. \qed

\end{document}